\documentclass[a4paper]{amsart}
\usepackage[enableskew,vcentermath]{youngtab}
\usepackage[dvips,line,matrix,frame,curve,poly]{xy}
\usepackage{amsthm,amsmath,amssymb,xspace,url}
\newtheorem{thm}{Theorem}[section]
\newtheorem{lem}[thm]{Lemma}

\newtheorem{prop}[thm]{Proposition}

\theoremstyle{definition}
\newtheorem{dfn}[thm]{Definition}

\theoremstyle{remark}
\newtheorem*{rmk}{Remark}

\newtheorem*{eg*}{Example}

\newcommand{\Set}[1]{\ensuremath{\mathcal{#1}}}            
\newcommand{\Mat}[1]{\ensuremath{\mathbf{#1}}}             
\newcommand{\Dfn}[1]{\emph{#1}}                            

\let\size=\abs

\newcommand{\jdt}{\emph{jeu de taquin}\xspace}
\newcommand{\J}{\overline j}

\def\alphaP{\alpha'}\def\betaP{\beta'}\def\gammaP{\gamma'}\def\deltaP{\delta'}
\def\alphaB{\bar\alpha}\def\betaB{\bar\beta}\def\gammaB{\bar\gamma}\def\deltaB{\bar\delta}
\def\alphaBP{\bar\alpha'}\def\betaBP{\bar\beta'}\def\gammaBP{\bar\gamma'}\def\deltaBP{\bar\delta'}

\author{Martin Rubey}
\address{Institut f\"ur Algebra, Zahlentheorie und Diskrete Mathematik, Leibniz
  Universit\"at Hannover, Welfengarten 1, D-30167 Hannover, Germany}
\email{martin.rubey@math.uni-hannover.de}
\urladdr{http://www.iazd.uni-hannover.de/~rubey/}
\thanks{Research partially supported by the Austrian Science Foundation FWF,
  grant S9607-N13, in the framework of the National Research Network \lq\lq
  Analytic Combinatorics and Probabilistic Number Theory\rq\rq.}


\title{Increasing and Decreasing Sequences in Fillings of Moon Polyominoes}

\usepackage[breaklinks=true,hyperref]{hyperref}
\begin{document}
\begin{abstract}
  We present an adaptation of \jdt and promotion for arbitrary fillings of moon
  polyominoes.  Using this construction we show various symmetry properties of
  such fillings taking into account the lengths of longest increasing and
  decreasing chains.  In particular, we prove a conjecture of Jakob Jonsson. We
  also relate our construction to the one recently employed by Christian
  Krattenthaler, thus generalising his results.
\end{abstract}
\maketitle

\def\dr{\POS[];[d]**\dir{-},[r]**\dir{-}}
\def\r{\POS[];[r]**\dir{-}}
\def\d{\POS[];[d]**\dir{-}}
\def\dd{\POS[];[d]**\dir{.}}
\def\ddr{\POS[];[d]**\dir{.},[r]**\dir{-}}
\def\e{\emptyset}
\def\x{\mbox{\Huge$\times$}}
\def\y{\POS[];[dr]**{}?*{\x}}
\def\c{\POS[];[dr]**{}?*{\mbox{\Huge$\circ$}}}
\def\put#1{\POS[];[dr]**{}?*{\mbox{#1}}}
\def\ann#1{\POS[];[dl]**{}?(0.25)*{#1}}
\def\B#1{\POS[];[dl]**{}?(0.25)*{\text{B#1}}}
\def\fB#1#2{\POS[];[dl]**{}?(0.25)*{\frac{\text{B#1}}{\text{B#2}}}}

\section{Introduction}
This article exploits properties of jeu de taquin, promotion and evacuation,
extended to fillings of matrices with non-negative integer entries.  Our
principal motivation was to prove a conjecture due to Jakob
Jonsson~\cite{Jonsson2005}, which is Theorem~\ref{thm:Jakob} of this article.
In rough terms, this theorem states that the number of fillings with zeros and
ones of a given moon polyomino (see Definition~\ref{dfn:moon}) that do not
contain a north-east chain of non-zero entries longer than a given threshold
(see Definition~\ref{sec:fillings-chains}) only depends on the distribution of
heights of the columns of the polyomino.  Aesthetics aside, one of the reasons
to consider this theorem is that for polyominoes of certain shapes the number
of such fillings is much easier to count than for others.

Jakob Jonsson's starting point was the desire to count generalised
triangulations with a given size of a maximal crossing.  He related these
objects to fans of Dyck paths, which were already counted by Myriam
de~Sainte-Catherine and Xavier Viennot~\cite{DeSainteCatherineViennot1986}.

As Christian Krattenthaler~\cite{Krattenthaler2006} noticed, Jakob Jonsson's
observation can be seen as a generalisation of symmetry properties of sizes of
maximal crossings and nestings in matchings and set-partitions.  These have
been proved a little earlier by William Chen, Eva Deng, Rosena Du, Richard
Stanley and Catherine Yan~\cite{ChenDengDuStanleyYan2006}, using the
Robinson-Schensted algorithm.  In my opinion, Christian Krattenthaler's most
important contribution was to prove these statements -- and generalisations
thereof -- using Sergey Fomin's growth
diagrams~\cite{Fomin1986,Fomin1994,Fomin1995,Roby1991}.

Curiously, all we use in this article are well known tools, although it seems
that some of their properties we need went unnoticed so far.  In particular,
Proposition~\ref{prop:equivalence} states a locality property of promotion and
the Robinson-Schensted-Knuth algorithm that might be interesting in its own
right.  A weaker form of this property has been used by Astrid
Reifegerste~\cite{Reifegerste2004} to deduce the sign of a permutation directly
from the pair of tableaux associated to it via the Robinson-Schensted
correspondence.

Meanwhile at least two articles dealing with other properties of fillings of
moon polyominoes have appeared.  In particular, Anisse
Kasraoui~\cite{Kasraoui2009} found that the \emph{number} of north-east and
south-east chains of length two in $0$-$1$-fillings of a moon polyomino remains
invariant if the columns of the polyomino are permuted, given that the result
is again a moon polyomino, and there is at most one non-zero entry in every
column.  Very remarkably, Anisse Kasraoui even found an astonishingly simple
expression for the generating function counting the number of such fillings.

William Chen, Svetlana Poznanovi\'c, Catherine Yan and Arthur
Yang~\cite{ChenPoznanovicYanYang2009} defined a major index for fillings of
$0$-$1$-fillings, that specialises to the major index for words and
permutations when the polyomino is a rectangle, and to the major index for
matchings and set partitions when the polyomino is a Ferrers shape.  Again,
this major index is invariant under permutations of the columns of the moon
polyomino, when the result is a moon polyomino.  Moreover, they found that the
generating function for their major index coincides with the one given by
Anisse Kasraoui.

In this article, we used the growth diagram description of the
Robinson-Schensted algorithm, and its Greene invariant that encodes in
particular the length of the longest north-east chain in the filling.  At least
for Theorem~\ref{thm:JakobWeak}, that deals with arbitrary fillings, an
analogous theorem should hold for any Greene invariant of a pair of dual-graded
graphs, where a promotion operator is available.  For example, we could take as
pair of dual graded graphs the lifted binary tree and the graph associated to
binary words, as in Section~4.6 of Sergey Fomin's article~\cite{Fomin1995}.
The resulting statistic on the fillings is then the maximal number of descents
(i.e., north-west chains of length two in consecutive rows) in a rectangle.  Of
course, the details have yet to be worked out, but there is certainly a lot
more to come.

\subsection*{Acknowledgements}
This article has undergone substantial changes since its first version, and I
would like to thank those who made this final version possible.  First of all
however I need to thank Christian Krattenthaler for presenting the problem in
the Arbeitsgemeinschaft Diskrete Mathematik, and then my wife Anita for
creating an atmosphere that let me have the crucial ideas during January 2006.
I feel indebted to several anonymous referees who made me check all my proofs
and thus discover (and fix!) many little and not so little gaps.  Last but not
least, I'm very grateful for the patience of \'editrice Mireille
Bousquet-M\'elou, and for constant support and interest of many colleagues.

\section{Definitions}

\subsection{Polyominoes}
\label{sec:polyominoes}
\begin{figure}[h]
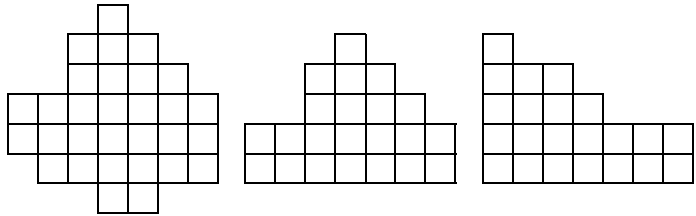

  \begin{equation*}
  \begin{array}{ccc}
  \young(:::\hfil,%
         ::\hfil\hfil\hfil,%
         ::\hfil\hfil\hfil\hfil,%
         \hfil\hfil\hfil\hfil\hfil\hfil\hfil,%
         \hfil\hfil\hfil\hfil\hfil\hfil\hfil,%
         :\hfil\hfil\hfil\hfil\hfil\hfil,%
         :::\hfil\hfil)  
  &
  \young(:::\hfil,%
         ::\hfil\hfil\hfil,%
         ::\hfil\hfil\hfil\hfil,%
         \hfil\hfil\hfil\hfil\hfil\hfil\hfil,%
         \hfil\hfil\hfil\hfil\hfil\hfil\hfil) 
  &
  \young(\hfil,%
         \hfil\hfil\hfil,%
         \hfil\hfil\hfil\hfil,%
         \hfil\hfil\hfil\hfil\hfil\hfil\hfil,%
         \hfil\hfil\hfil\hfil\hfil\hfil\hfil)
  \end{array}
  \end{equation*}
  \caption{a moon-polyomino, a stack-polyomino and a Ferrers diagram}
  \label{fig:moon}
\end{figure}

\begin{dfn}\label{dfn:polyominoes}
  A \Dfn{polyomino} is a finite subset of $\mathbb Z^2$, where we regard an
  element of $\mathbb Z^2$ as a cell.  A \Dfn{column} of a polyomino is the set
  of cells along a vertical line, a \Dfn{row} is the set of cells along a
  horizontal line.

  The polyomino is \Dfn{convex}, if for any two cells in a column, the elements
  of $\mathbb Z^2$ in between are also cells of the polyomino, and for any two
  cells in a row, the elements of $\mathbb Z^2$ in between are also cells of
  the polyomino.  It is \Dfn{intersection-free}, if any two columns are
  \Dfn{comparable}, i.e., the set of row coordinates of cells in one column is
  contained in the set of row coordinates of cells in the other.  Equivalently,
  it is intersection-free, if any two rows are comparable.

  For example, the polyomino
  \begin{equation*}
    \young(::\hfil,%
    ::\hfil\hfil\hfil,%
    \hfil\hfil\hfil\hfil\hfil,%
    \hfil\hfil\hfil\hfil,%
    ::\hfil)  
  \end{equation*}
  is convex, but not intersection-free, since the first and the last columns are
  incomparable.
\end{dfn}

\begin{dfn}\label{dfn:moon}
  A \Dfn{moon polyomino} is a convex, intersection-free polyomino.  A
  \Dfn{stack polyomino} is a moon-polyomino if all columns start at the same
  level.  A \Dfn{Ferrers diagram} is a stack-polyomino with weakly decreasing
  row widths $\lambda_1,\lambda_2,\dots,\lambda_n$, reading rows from bottom to
  top.  We alert the reader that we are using \lq French\rq\ notation for
  Ferrers diagrams.

  For a permutation $\sigma$ of the column indices of a moon polyomino $M$,
  \Dfn{reordering the columns according to $\sigma$} yields the polyomino
  $\sigma M=\{(\sigma i,j): (i,j)\in M\}$.  Thus, the columns are only
  translated horizontally and keep their vertical offsets.

  The \Dfn{content} of a moon polyomino is the sequence of column heights, in
  decreasing order.  For example, the content of the moon polyomino at the left
  of Figure~\ref{fig:moon} is $(7,6,5,4,3,3,2)$, while the content of the other
  two polyominoes in the same figure is $(5,4,4,3,2,2,2)$.
\end{dfn}
\begin{rmk}
  The name \lq moon polyomino\rq\ was used by Jakob Jonsson
  in~\cite{Jonsson2005}.  Curiously, this class of polyominoes was
  independently introduced a little earlier by Giusi Castiglione and Antonio
  Restivo~\cite{CastiglioneRestivo2003} as \Dfn{L-convex} (or, alternatively
  \Dfn{1-convex}) polyominoes.  The defining property they use is that every
  pair of cells can be connected by a path consisting of horizontal and
  vertical steps that changes direction at most once.
\end{rmk}

\subsection{Fillings and Chains}
\label{sec:fillings-chains}
\begin{dfn}\label{dfn:fillings}
  An \Dfn{arbitrary filling} of a polyomino is an assignment of natural numbers
  to the cells of the polyomino.  We refer to the number in a cell as the
  \Dfn{multiplicity} of the entry.

  In a \Dfn{$0$-$1$-filling} we restrict ourselves to the numbers $0$ and $1$.
  A \Dfn{standard filling} has the additional constraint that in each column
  and in each row there is exactly one entry $1$, whereas a \Dfn{partial
    filling} has at most one entry $1$ in each column and in each row.
\end{dfn}
\begin{rmk}
  In the figures, we will usually omit zeros, and in $0$-$1$-fillings we will
  replace ones by crosses for \ae sthetic reasons.  We reserve the letters
  $\alpha$, $\beta$, $\gamma$, $\delta$, $\epsilon$ and $\pi$ to denote
  fillings.
\end{rmk}

\begin{dfn}
  A \Dfn{chain} is a sequence of non-zero entries in a filling such that the
  smallest rectangle containing all its elements of the sequence is completely
  contained in the moon polyomino.

  A \Dfn{north-east chain}, or short \Dfn{ne-chain} of length $k$ in an
  arbitrary filling of a moon polyomino is a chain of $k$ non-zero entries,
  such that each element is strictly to the right and strictly above the
  preceding element of the sequence.  Similarly, in a \Dfn{south-east chain},
  for short \Dfn{se-chain}, each element is strictly to the right and strictly
  below the preceding element.  The length of such a chain is the number of its
  elements.

  \Dfn{NE-chains} and \Dfn{SE-chains} may have elements in the same column and
  in the same row.  For these kinds of chains, each element contributes its
  multiplicity to the length.  That is, a \Dfn{NE-chain} of length $k$ is a
  chain such that each element is weakly to the right and weakly above its
  predecessor, and the sum of the multiplicities of the elements equals $k$.

  For $0$-$1$-fillings we also define \Dfn{nE-chains} and \Dfn{sE-chains},
  where we allow an element of the chain to be in the same column as its
  predecessor, but not in the same row.  Similarly, elements of \Dfn{Ne-chains}
  and \Dfn{Se-chains} are allowed to be in the same row, but not in the same
  column.
\end{dfn}
For example, consider the following two fillings:
\begin{equation*}
\young(:\hfil1,%
       :\hfil13,%
       3\hfil1\hfil,%
       11\hfil\hfil)\quad\text{and}\quad
\young(:\hfil1,%
       :\hfil1\hfil,%
       \hfil\hfil1\hfil,%
       \hfil1\hfil\hfil)
\end{equation*}
The length of the longest ne-chain in the filling on the left is three, whereas
the length of the longest se-chain is two. The lengths of the longest NE- and
SE-chains are six and five respectively.

The lengths of the longest nE-, Ne-, sE-, and Se-chains in the $0$-$1$-filling
on the right are four, two, three and one respectively.

\subsection{Partitions and Tableaux}
\label{sec:partitions-tableaux}
\begin{dfn}
  A \Dfn{partition} is a weakly decreasing sequence of natural numbers, which
  are called its \Dfn{parts}. The \Dfn{length} of a partition is the number of
  its parts, the \Dfn{size} of a partition is the sum of its parts. A partition
  $\lambda=(\lambda_1,\lambda_2,\dots,\lambda_l)$ is \Dfn{contained} in another
  partition $\mu=(\mu_1,\mu_2,\dots,\mu_m)$ if $l\leq m$ and
  $\lambda_i\leq\mu_i$ for all $i\leq l$.

  The union of $\lambda$ and $\mu$, denoted $\lambda\cup\mu$, is the partition
  $\kappa=(\kappa_1,\kappa_2,\dots,\kappa_k)$ with $k=\max(l,m)$ and
  $\kappa_i=\max(\lambda_i, \mu_i)$ for all $i\leq k$, where we set
  $\lambda_i=0$ for $i>l$ and $\mu_i=0$ for $i>m$.

  The \Dfn{transpose} or \Dfn{conjugate} of a partition $\lambda$ is defined as
  $\lambda^t=(\mu_1,\mu_2,\dots,\mu_m)$, where $m=\lambda_1$ and $\mu_i$ is the
  number of parts in $\lambda$ greater than or equal to $i$.

  The \Dfn{transpose} of a sequence of partitions
  $P=(\emptyset=\lambda^0,\lambda^1,\dots,\lambda^n)$ is the sequence of
  partitions $P^t$ obtained by transposing each individual partition.
\end{dfn}

\begin{rmk}
  A Ferrers shape (in French notation) corresponds to a partition $\lambda$,
  setting $\lambda_i$ to the length of the $i$\textsuperscript{th} row from
  bottom to top.  Using this correspondence, the transpose of a partition can
  be obtained by reflecting the corresponding Ferrers shape about the main
  diagonal.
\end{rmk}

\begin{dfn}
  A \Dfn{semi-standard Young tableau} is a filling of a Ferrers shape with
  positive integers, such that entries are weakly increasing in rows and
  strictly increasing -- from bottom to top -- in columns.
 
  A \Dfn{standard Young tableau} is a semi-standard Young tableau with entries
  being the numbers $1$ through $n$, such that each number occurs exactly once.

  A \Dfn{partial Young tableau} is a semi-standard Young tableau with all
  entries distinct.
\end{dfn}
\begin{rmk}
  An example for a pair of standard Young tableaux is given in
  Equation~\eqref{eq:Young-tableaux}.  Although Young tableaux happen to be
  fillings of Ferrers shapes, they should be thought of as objects entirely
  different from the fillings introduced in Definition~\ref{dfn:fillings}.
\end{rmk}

\section{Growth Diagrams and the Robinson-Schensted-Knuth
  Correspondence}\label{sec:growth}

Sergey Fomin's growth diagrams together with Marcel Sch\"utzenberger's \jdt and
\emph{promotion} will be the central tools in this article.  In this section we
introduce growth diagrams, which associate sequences of partitions to fillings
of matrices with non-negative integer entries.  Although the contents of this
section is well known, we reproduce it here for the convenience of the
reader. Some additional information and more references can be found in
\cite[Sections~2 and 4]{Krattenthaler2006}.

\subsection{Local Rules and Growth Diagrams}\label{sec:local-rules}
Consider a rectangular polyomino with a partial filling, as, for example, in
Figure~\ref{fig:growth}.a where we have replaced zeros by empty cells and ones
by crosses.  Using the following construction we will inductively label the
corners of each cell with a partition, starting from the bottom left corner, to
obtain a \Dfn{growth diagram}.

\begin{figure}[htb]
  \begin{equation*}
\begin{xy}
  \xymatrix@!{%
    \mu    \dr&\rho\d\\
    \lambda\r &\nu}
\end{xy}
  \end{equation*}
  \caption{a cell of a growth diagram}
  \label{fig:growth-cell}
\end{figure}
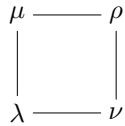

First, we attach the empty partition $\emptyset$ to the corners along the lower
and the left border.  Suppose now that we have already labelled all the corners
of a square except the top right with partitions $\lambda$, $\mu$ and $\nu$, as
in Figure~\ref{fig:growth-cell}. We compute $\rho$ as follows:

\begin{itemize}
\item[F1] Suppose that the square does not contain a cross, and that
  $\lambda=\mu=\nu$. Then set $\rho=\lambda$.
\item[F2] Suppose that the square does not contain a cross, and that
  $\mu\neq\nu$. Then set $\rho=\mu\cup\nu$.
\item[F3] Suppose that the square does not contain a cross, and that
  $\lambda\subsetneq\mu=\nu$. Then we obtain $\rho$ from $\mu$ by adding $1$ to
  the $i+1$\textsuperscript{st} part of $\mu$, given that $\lambda$ and $\mu$
  differ in the $i$\textsuperscript{th} part.
\item[F4] Suppose that the square contains a cross. This implies that
  $\lambda=\mu=\nu$ and we obtain $\rho$ from $\lambda$ by adding $1$ to the
  first part of $\lambda$.
\end{itemize}

\begin{figure}[h]
\begin{equation*}
\begin{xy}
  \POS(0,-42)
  \drop{a.}
  \POS(12,0)
  \xymatrix@!=4pt{%
  \e\dr  &1 \dr  &11\dr  &21\dr  &211\dr  &311\dr\y&411\dr  &421\dr  &4211\d\\
  \e\dr  &1 \dr  &11\dr  &21\dr  &211\dr\y&311\dr  &311\dr  &321\dr  &3211\d\\
  \e\dr\y&1 \dr  &11\dr  &21\dr  &211\dr  &211\dr  &211\dr  &311\dr  &3111\d\\
  \e\dr  &\e\dr  &1 \dr\y&2 \dr  &21 \dr  &21 \dr  &21 \dr  &31 \dr  &311 \d\\
  \e\dr  &\e\dr  &1 \dr  &1 \dr  &2  \dr  &2  \dr  &2  \dr\y&3  \dr  &31  \d\\
  \e\dr  &\e\dr  &1 \dr  &1 \dr\y&2  \dr  &2  \dr  &2  \dr  &2  \dr  &21  \d\\
  \e\dr  &\e\dr  &1 \dr  &1 \dr  &1  \dr  &1  \dr  &1  \dr  &1  \dr\y&2   \d\\
  \e\dr  &\e\dr\y&1 \dr  &1 \dr  &1  \dr  &1  \dr  &1  \dr  &1  \dr  &1   \d\\
  \e\r   &\e\r   &\e\r   &\e\r   &\e \r   &\e \r   &\e \r   &\e \r   &\e}
  \POS(0,-142)
  \drop{b.}
  \POS(12,-100)
  \xymatrix@!=4pt{%
&\e\dr  &1 \dr  &2 \dr  &21 \dr  &31 \dr\y&41 \dr  &42  \dr  &421 \dr  &4211\d\\
&{}\dr  &{}\dr  &{}\dr  &{} \dr\y&{} \dr  &{} \dr  &{}  \dr  &{}  \dr  &3211\d\\
&{}\dr  &{}\dr\y&{}\dr  &{} \dr  &{} \dr  &{} \dr  &{}  \dr  &{}  \dr  &3111\d\\
&{}\dr  &{}\dr  &{}\dr  &{} \dr  &{} \dr  &{} \dr\y&{}  \dr  &{}  \dr  &311 \d\\
&{}\dr  &{}\dr  &{}\dr  &{} \dr  &{} \dr  &{} \dr  &{}  \dr\y&{}  \dr  &31  \d\\
&{}\dr  &{}\dr  &{}\dr\y&{} \dr  &{} \dr  &{} \dr  &{}  \dr  &{}  \dr  &21  \d\\
&{}\dr  &{}\dr  &{}\dr  &{} \dr  &{} \dr  &{} \dr  &{}  \dr  &{}  \dr\y&2   \d\\
&{}\dr\y&{}\dr  &{}\dr  &{} \dr  &{} \dr  &{} \dr  &{}  \dr  &{}  \dr  &1   \d\\
&{}\r   &{}\r   &{}\r   &{} \r   &{} \r   &{} \r   &{}  \r   &{}  \r   &\e}
\end{xy}
\end{equation*}
\caption{a. a growth diagram\quad b. \emph{promotion} on the upper border}
\label{fig:growth}
\end{figure}

The important fact is, that this process is invertible: given the labels of the
corners along the upper and right border of the diagram, we can reconstruct the
complete growth diagram as well as the entries of the squares.  To this end,
suppose that we have already labelled all the corners of a square except the
bottom left with partitions $\mu$ and $\nu$ and $\rho$, as in
Figure~\ref{fig:growth-cell}. We compute $\lambda$ and the entry of the square
as follows:

\begin{itemize}
\item[B1] If $\mu=\nu=\rho$ we set $\lambda=\rho$ and leave the square empty.
\item[B2] If $\mu\neq\nu$ we set $\lambda=\mu\cap\nu$ and leave the square empty.
\item[B3] If $\mu=\nu\subsetneq\rho$ and $\mu$ and $\rho$ differ in the
  $i$\textsuperscript{th} part for $i\geq2$, we obtain $\lambda$ from $\mu$ by
  deleting $1$ from the $i-1$\textsuperscript{st} part of $\mu$ and leave the
  square empty.
\item[B4] If $\mu=\nu\subsetneq\rho$ and $\mu$ and $\rho$ differ in the first
  part we set $\lambda=\mu$ and mark the square with a cross.
\end{itemize}

\subsection{The Robinson-Schensted Correspondence and Greene's
  Theorem}\label{sec:RSK-Greene}
In the case of a standard filling of a square, the sequence of partitions
$\emptyset=\mu^0,\mu^1,\dots,\mu^n$ along the upper border of the growth
diagram corresponds to a standard Young tableau $Q$ as follows: we put the
entry $i$ into the cell by which $\mu^{i-1}$ and $\mu^i$ differ.  Similarly,
the sequence of partitions $\emptyset=\lambda^0,\lambda^1,\dots,\lambda^n$
along the right border of the diagram corresponds to a standard Young tableau
$P$ of the same shape as $Q$.
  
Furthermore, the filling itself defines a permutation $\pi$. For example in
Figure~\ref{fig:growth}.a we have $\pi=6,1,5,3,7,8,4,2$ and
\begin{equation}\label{eq:Young-tableaux}
  \Yinterspace0.6ex plus 0.3ex
  (P,Q)=\left(\young(6,5,37,1248),\young(8,4,27,1356)\right).
\end{equation}
It is well known (see for example Theorem 7.13.5 of \cite{EC2}) that $Q$ is the
recording and $P$ the insertion tableau produced by the Robinson-Schensted
correspondence, applied to the permutation $\pi$.

Since the partitions along the upper and right border of a growth diagram
determine the filling and vice versa, the following definition will be useful:
\begin{dfn}
  Let $\pi$ be a partial filling of a rectangular polyomino and consider the
  corresponding growth diagram.  Suppose that the corners along the right
  border are labelled with a sequence of partitions $P$, and along the upper
  border with a sequence of partitions $Q$.  We then say, that $\pi$
  \Dfn{corresponds} to the pair $(P, Q)$.  The last partitions of $P$ and $Q$
  coincide, and this partition is the \Dfn{shape} of the filling $\pi$.
\end{dfn}

For our purposes it is of great importance that the partitions appearing in the
corners of a growth diagram also have a \lq global\rq\ description. This is
called \Dfn{Greene's Theorem} (due to Curtis Greene~\cite{Greene1974}):
\begin{thm}\cite[Theorem A.1.1.1]{EC2}\label{thm:Greene}
  Suppose that a corner $c$ of a growth diagram is labelled by a partition
  $\lambda$.  Then, for any integer $k$, the maximal cardinality of a union of
  $k$ pairwise disjoint north-east chains situated in the rectangular region to
  the left and below of $c$ is equal to $\lambda_1+\lambda_2+\dots+\lambda_k$.

  Furthermore, the maximal cardinality of the union of $k$ pairwise disjoint
  south-east chains situated in the rectangular region to the left and below of
  $c$ is equal to $\lambda'_1+\lambda'_2+\dots+\lambda'_k$, where $\lambda'$ is
  the transpose of $\lambda$.
\end{thm}
\begin{rmk}
  In \cite{EC2} this theorem is only stated for standard fillings.  To obtain
  the statement for partial fillings, it suffices to observe that by rule F1
  empty columns and rows can be ignored.
\end{rmk}

\subsection{Variations of the Robinson-Schensted-Knuth Correspondence}
\label{sec:RS-variations}
In the following we extend the construction described in
Section~\ref{sec:local-rules} to arbitrary fillings of rectangular polyominoes.
We construct a new diagram with more rows and columns, and place non-zero
entries which are originally in the same column or row in different columns and
rows in the larger diagrams.  A similar strategy is applied to entries with
multiplicity larger than one.  We refer to this process, transforming a filling
$\pi$ into a standard filling $std(\pi)$, as \Dfn{standardisation}.  More
precisely, we proceed as follows:

\subsubsection{First Variant: RSK}
\label{sec:first-variant:-rsk}
\begin{figure}[h]
  \begin{equation*}
  \begin{xy}<0pt,-50pt>;<10pt,-50pt>:
    (0,0);(4,0)**[|(3)]\dir{-},
    (0,2);(4,2)**[|(3)]\dir{-},
    (0,4);(4,4)**[|(3)]\dir{-},
    (0,6);(4,6)**[|(3)]\dir{-},
    (0,8);(4,8)**[|(3)]\dir{-},
    (0,0);(0,8)**[|(3)]\dir{-},
    (2,0);(2,8)**[|(3)]\dir{-},
    (4,0);(4,8)**[|(3)]\dir{-},
    (1,1),*{3},(3,5),*{1},(1,7)*{1},(3,7)*{1}
  \end{xy}
  \quad\mapsto\quad
  \begin{xy}<0pt,-70pt>;<10pt,-70pt>:
    (0,0);(12,0)**[|(3)]\dir{-},
    (0,2);(12,2)**\dir{-},
    (0,4);(12,4)**\dir{-},
    (0,6);(12,6)**[|(3)]\dir{-},
    (0,8);(12,8)**[|(3)]\dir{-},
    (0,10);(12,10)**\dir{-},
    (0,12);(12,12)**[|(3)]\dir{-},
    (0,0);(0,12)**[|(3)]\dir{-},
    (2,0);(2,12)**\dir{-},
    (4,0);(4,12)**\dir{-},
    (6,0);(6,12)**\dir{-},
    (8,0);(8,12)**[|(3)]\dir{-},
    (10,0);(10,12)**\dir{-},
    (12,0);(12,12)**[|(3)]\dir{-},
    (7,9),*{\times},(11,11),*{\times},(9,7),*{\times},
    (5,5),*{\times},(3,3),*{\times},(1,1),*{\times},
    (0.4,12.4),*[r]{\scriptstyle\emptyset},
    (8.4,12.4),*[r]{\scriptstyle 4},
   (12.4,12.4),*[r]{\scriptstyle 51},
    (8.4,12.4),*[r]{\scriptstyle 4},
   (12.4,8.4),*[r]{\scriptstyle 4},
   (12.4,6.4),*[r]{\scriptstyle 3},
   (12.4,0.4),*[r]{\scriptstyle\emptyset}
  \end{xy}
  \;\text{or}\quad
  \begin{xy}<0pt,-70pt>;<10pt,-70pt>:
    (0,0);(12,0)**[|(3)]\dir{-},
    (0,2);(12,2)**\dir{-},
    (0,4);(12,4)**\dir{-},
    (0,6);(12,6)**[|(3)]\dir{-},
    (0,8);(12,8)**[|(3)]\dir{-},
    (0,10);(12,10)**\dir{-},
    (0,12);(12,12)**[|(3)]\dir{-},
    (0,0);(0,12)**[|(3)]\dir{-},
    (2,0);(2,12)**\dir{-},
    (4,0);(4,12)**\dir{-},
    (6,0);(6,12)**\dir{-},
    (8,0);(8,12)**[|(3)]\dir{-},
    (10,0);(10,12)**\dir{-},
    (12,0);(12,12)**[|(3)]\dir{-},
    (1,11),*{\times},(9,9),*{\times},(11,7),*{\times},
    (3,5),*{\times},(5,3),*{\times},(7,1),*{\times},
    (0.4,12.4),*[r]{\scriptstyle\emptyset},
    (8.4,12.4),*[r]{\scriptstyle 1111},
    (12.4,12.4),*[r]{\scriptstyle 2211},
    (12.4,8.4),*[r]{\scriptstyle 211},
    (12.4,6.4),*[r]{\scriptstyle 111},
    (12.4,0.4),*[r]{\scriptstyle\emptyset}
  \end{xy}
  \end{equation*}
  \caption{standardisation of an arbitrary filling using RSK or dual RSK'}
  \label{fig:blow-up-RSK-Burge}
\end{figure}
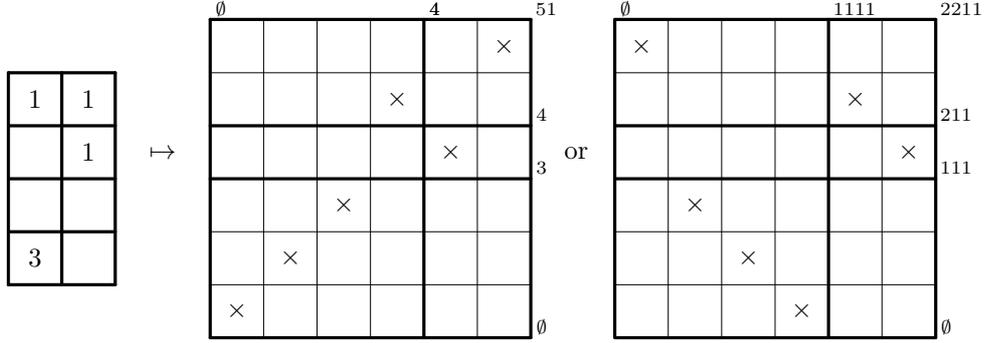

Each row and each column of the original diagram is replaced by as many rows
and columns in the new diagram as it contains non-zero entries, counting
multiplicities.  Then, for each row and for each column of the original diagram
we place the non-zero entries into the new diagram as a north-east chain.  An
example can be found in Figure~\ref{fig:blow-up-RSK-Burge}, the result being
the left of the two standardised diagrams.  Note that this process preserves
the length of the NE- and se-chains.

Given a filling $\pi$, we can apply the rules F1 to F4 to the transformed
diagram and obtain a pair of sequences of partitions $(\bar P, \bar Q)$ with
$\bar P=(\emptyset=\bar\lambda^0, \bar\lambda^1, \dots, \bar\lambda^n)$ and
$\bar Q=(\emptyset=\bar\mu^0, \bar\mu^1, \dots, \bar\mu^n)$.  The sequences of
partitions corresponding to the original diagram are then given by
$P=(\emptyset=\lambda^0, \lambda^1, \dots, \lambda^k)$ and $Q=(\emptyset=\mu^0,
\mu^1, \dots, \mu^k)$ with $\lambda^i = \bar\lambda^{\#\text{non-zero entries
    below row $i$}}$ and $\mu^i = \bar\mu^{\#\text{non-zero entries left of
    column $i$}}$.  Less formally: we keep only the partitions between \lq
compartments\rq, and partitions are constant across empty rows and columns.

We remark that the partitions $\lambda^{i-1}$ and $\lambda^i$, as well as
$\mu^{i-1}$ and $\mu^i$ differ by a \Dfn{horizontal strip} each, i.e., by at
most one cell in each column.  Thus, $P$ corresponds to a semi-standard Young
tableau as follows: we put an entry $i$ into each cell that is present in
$\lambda^i$ but not in $\lambda^{i-1}$.  Similarly, also $Q$ corresponds to a
semi-standard Young tableau.  In the example of
Figure~\ref{fig:blow-up-RSK-Burge} (left diagram), we have
\begin{equation*}
  \Yinterspace0.6ex plus 0.3ex
  (P,Q)=\left(\young(3,11123),\young(2,11112)\right).  
\end{equation*}
It is well known that $(P,Q)$ coincides with the result of applying the usual
\lq Robinson-Schensted-Knuth\rq, short RSK correspondence, to $\pi$.

We stress that the partitions in $\bar P$ and $\bar Q$ can be reconstructed
from those in $P$ and $Q$ alone, given the knowledge that the non-zero entries
of the filling within compartments should be arranged in north-east chains.  In
fact, the local rules F1 to F4 and B1 to B4 can be adapted to work
\emph{directly} on the original diagrams, bypassing the standardisation process
entirely.  The details can be found in Christian Krattenthaler's article
\cite{Krattenthaler2006}, Section~4.  For convenience, we restate the
corresponding variant of Greene's theorem from \cite{Krattenthaler2006},
Section~4:
\begin{thm}\label{thm:Greene-RSK}
  Suppose that a corner $c$ of a growth diagram completed according to the RSK
  correspondence is labelled by a partition $\lambda$.  Then, for any integer
  $k$, the maximal cardinality of a union of $k$ pairwise disjoint NE-chains
  situated in the rectangular region to the left and below of $c$ is equal to
  $\lambda_1+\lambda_2+\dots+\lambda_k$.

  Furthermore, the maximal cardinality of the multiset union of $k$ se-chains,
  where any entry $e$ of the filling appears in at most $e$ chains, all
  situated in the rectangular region to the left and below of $c$, is equal to
  $\lambda'_1+\lambda'_2+\dots+\lambda'_k$.
\end{thm}

\subsubsection{Second Variant: dual RSK'}
\label{sec:second-variant:-dual}
There is another obvious possibility to standardise an arbitrary filling.
Instead of placing the non-zero entries into the new diagram as north-east
chains, we could also arrange them in south-east chains, thus preserving the
length of ne- and SE-chains.  An example for this transformation is given in
Figure~\ref{fig:blow-up-RSK-Burge}, the result being the right of the two
standardised diagrams.

In this case, the partitions $\lambda^{i-1}$ and $\lambda^i$ as well as
$\mu^{i-1}$ and $\mu^i$ differ by a \Dfn{vertical strip} each, i.e., by at most
one cell in each row.  Putting $i$ into each cell that is present in $\lambda^i$
but not in $\lambda^{i-1}$, we see that $P$ corresponds to a \emph{transposed}
semi-standard Young tableau, similarly for $Q$.  In the example of
Figure~\ref{fig:blow-up-RSK-Burge} (right diagram), we obtain
\begin{equation*}
  \Yinterspace0.6ex plus 0.3ex
  (P,Q)=\left(\young(3,1,13,12),\young(1,1,12,12)\right).
\end{equation*}
What we just described is the growth diagram version of the dual RSK'
correspondence, also known as the \lq Burge\rq\ correspondence.

Greene's theorem now reads as follows:
\begin{thm}\label{thm:Greene-Burge}
  Suppose that a corner $c$ of a growth diagram completed according to the dual
  RSK' correspondence is labelled by a partition $\lambda$.  Then, for any
  integer $k$, the maximal cardinality of a multiset union of $k$ ne-chains
  where any entry $e$ of the filling appears in at most $e$ chains, all
  situated in the rectangular region to the left and below of $c$ is equal to
  $\lambda_1+\lambda_2+\dots+\lambda_k$.

  Furthermore, the maximal cardinality of the union of $k$ pairwise disjoint
  SE-chains situated in the rectangular region to the left and below of $c$, is
  equal to $\lambda'_1+\lambda'_2+\dots+\lambda'_k$.
\end{thm}

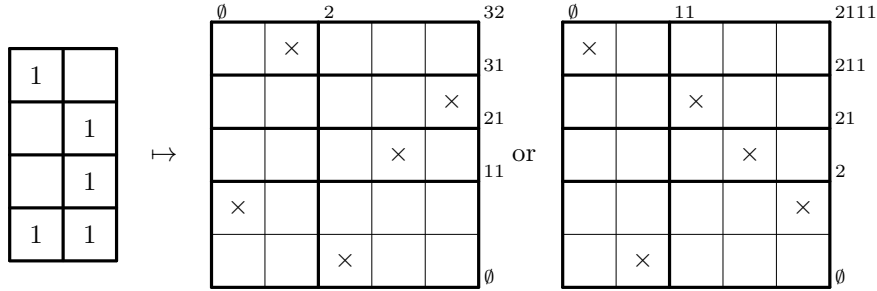
\begin{figure}[h]
  \begin{equation*}
  \begin{xy}<0pt,-40pt>;<10pt,-40pt>:
    (0,0);(4,0)**[|(3)]\dir{-},
    (0,2);(4,2)**[|(3)]\dir{-},
    (0,4);(4,4)**[|(3)]\dir{-},
    (0,6);(4,6)**[|(3)]\dir{-},
    (0,8);(4,8)**[|(3)]\dir{-},
    (0,0);(0,8)**[|(3)]\dir{-},
    (2,0);(2,8)**[|(3)]\dir{-},
    (4,0);(4,8)**[|(3)]\dir{-},
    (1,1),*{1},(1,7),*{1},(3,1),*{1},(3,3)*{1},(3,5)*{1}
  \end{xy}
  \quad\mapsto\quad
  \begin{xy}<0pt,-50pt>;<10pt,-50pt>:
    (0,0);(10,0)**[|(3)]\dir{-},
    (0,2);(10,2)**\dir{-},
    (0,4);(10,4)**[|(3)]\dir{-},
    (0,6);(10,6)**[|(3)]\dir{-},
    (0,8);(10,8)**[|(3)]\dir{-},
    (0,10);(10,10)**[|(3)]\dir{-},
    (0,0);(0,10)**[|(3)]\dir{-},
    (2,0);(2,10)**\dir{-},
    (4,0);(4,10)**[|(3)]\dir{-},
    (6,0);(6,10)**\dir{-},
    (8,0);(8,10)**\dir{-},
    (10,0);(10,10)**[|(3)]\dir{-},
    (1,3),*{\times},(3,9),*{\times},(5,1),*{\times},(7,5),*{\times},
    (9,7),*{\times},
    (0.4,10.4),*[r]{\scriptstyle\emptyset},
    (4.4,10.4),*[r]{\scriptstyle 2},
   (10.4,10.4),*[r]{\scriptstyle 32},
   (10.4,8.4),*[r]{\scriptstyle 31},
   (10.4,6.4),*[r]{\scriptstyle 21},
   (10.4,4.4),*[r]{\scriptstyle 11},
   (10.4,0.4),*[r]{\scriptstyle\emptyset}
  \end{xy}
  \;\text{or}\quad
  \begin{xy}<0pt,-50pt>;<10pt,-50pt>:
    (0,0);(10,0)**[|(3)]\dir{-},
    (0,2);(10,2)**\dir{-},
    (0,4);(10,4)**[|(3)]\dir{-},
    (0,6);(10,6)**[|(3)]\dir{-},
    (0,8);(10,8)**[|(3)]\dir{-},
    (0,10);(10,10)**[|(3)]\dir{-},
    (0,0);(0,10)**[|(3)]\dir{-},
    (2,0);(2,10)**\dir{-},
    (4,0);(4,10)**[|(3)]\dir{-},
    (6,0);(6,10)**\dir{-},
    (8,0);(8,10)**\dir{-},
    (10,0);(10,10)**[|(3)]\dir{-},
    (1,9),*{\times},(3,1),*{\times},(5,7),*{\times},(7,5),*{\times},
    (9,3),*{\times},
    (0.4,10.4),*[r]{\scriptstyle\emptyset},
    (4.4,10.4),*[r]{\scriptstyle 11},
   (10.4,10.4),*[r]{\scriptstyle 2111},
   (10.4,8.4),*[r]{\scriptstyle 211},
   (10.4,6.4),*[r]{\scriptstyle 21},
   (10.4,4.4),*[r]{\scriptstyle 2},
   (10.4,0.4),*[r]{\scriptstyle\emptyset}
  \end{xy}
  \end{equation*}
  \caption{standardisation of a $0$-$1$-filling using dual RSK or RSK'}
  \label{fig:blow-up-dual-RSK-RSK-prime}
\end{figure}

\subsubsection{Third Variant: dual RSK}
\label{sec:third-variant:-dual}
If we restrict ourselves to $0$-$1$-fillings, we can also transform multiple
non-zero entries of a column of the original diagram into a north-east chain
and multiple non-zero entries of a row into a south-east chain.  In this case,
the lengths of nE- and Se-chains are preserved, as can be seen from the example
on the left of Figure~\ref{fig:blow-up-dual-RSK-RSK-prime}.

Now $\lambda^{i-1}$ and $\lambda^i$ differs by a vertical strip, while
$\mu^{i-1}$ and $\mu^i$ differs by a horizontal strip.  We thus obtain the
so-called dual RSK correspondence.  In the example of
Figure~\ref{fig:blow-up-dual-RSK-RSK-prime} (left diagram), we obtain
\begin{equation*}
  \Yinterspace0.6ex plus 0.3ex
  (P,Q)=\left(\young(14,123),\young(22,112)\right),
\end{equation*} 
and $P$ is a transposed semi-standard Young tableau, while $Q$ is a
semi-standard Young tableau.

Again we note the corresponding variation of Greene's theorem:
\begin{thm}
  Suppose that a corner $c$ of a growth diagram completed according to the dual
  RSK correspondence is labelled by a partition $\lambda$.  Then, for any
  integer $k$, the maximal cardinality of a union of $k$ pairwise disjoint
  nE-chains situated in the rectangular region to the left and below of $c$ is
  equal to $\lambda_1+\lambda_2+\dots+\lambda_k$.

  Furthermore, the maximal cardinality of the union of $k$ pairwise disjoint
  Se-chains situated in the rectangular region to the left and below of $c$, is
  equal to $\lambda'_1+\lambda'_2+\dots+\lambda'_k$.
\end{thm}

\subsubsection{Fourth Variant: RSK'}
\label{sec:fourth-variant:-rsk}
As a last possibility, again for $0$-$1$-fillings, we can transform multiple
non-zero entries of a column of the original diagram into a south-east chain
and multiple non-zero entries of a row into a north-east chain, obtaining the
\lq Robinson-Schensted-Knuth-prime\rq\ correspondence, short RSK', which
preserves the lengths of Ne- and sE-chains. This is shown on the right of
Figure~\ref{fig:blow-up-dual-RSK-RSK-prime}.

Here $\lambda^{i-1}$ and $\lambda^i$ differs by a horizontal strip, while
$\mu^{i-1}$ and $\mu^i$ differs by a vertical strip.  In the example of
Figure~\ref{fig:blow-up-dual-RSK-RSK-prime} (right diagram), we obtain
\begin{equation*}
  \Yinterspace0.6ex plus 0.3ex
  (P,Q)=\left(\young(4,3,2,11),\young(2,2,1,12)\right),
\end{equation*} 
and $P$ is a semi-standard Young tableau, while $Q$ is a transposed
semi-standard Young tableau.

Thus, the last variation of Greene's theorem is:
\begin{thm}
  Suppose that a corner $c$ of a growth diagram completed according to the RSK'
  correspondence is labelled by a partition $\lambda$.  Then, for any integer
  $k$, the maximal cardinality of a union of $k$ pairwise disjoint Ne-chains
  situated in the rectangular region to the left and below of $c$ is equal to
  $\lambda_1+\lambda_2+\dots+\lambda_k$.

  Furthermore, the maximal cardinality of the union of $k$ pairwise disjoint
  sE-chains situated in the rectangular region to the left and below of $c$, is
  equal to $\lambda'_1+\lambda'_2+\dots+\lambda'_k$.
\end{thm}

\section{Jeu de Taquin and Promotion}
\label{sec:jeu-de-taquin}

Our second tool is \jdt, introduced by Marcel Sch\"utzenberger, an algorithm
that \lq rectifies\rq\ a skew semi-standard Young tableau.  For our purposes it
is not necessary to introduce skew tableaux, and what we call \jdt is sometimes
referred to as the \Dfn{$\Delta$-operator}, eg.\ in Bruce Sagan's
book~\cite{SymmetricGroup}.  Because of its simplicity, we begin with the
description for partial tableaux.

Given a partial Young tableau $P$, we define $jdt(P)$ as the result of the
following algorithm:
\begin{enumerate}
\item subtract one from all the entries in the tableau.
\item if present, replace the cell containing $0$ by an empty cell. Otherwise
  stop.  In the following, we will move the empty cell to the top right border
  of the tableau, given French notation.
\item if there is no cell to the right and no cell above the empty cell, remove
  the empty cell and stop.
\item consider the cells above and to the right of the cell without entry, and
  exchange the cell with the smaller entry and the empty cell. Go to Step~3.
\end{enumerate}

This variation of \emph{jeu de taquin} can be conveniently described with \lq
growth diagrams\rq\ -- of a different kind than those introduced in
Section~\ref{sec:local-rules} -- as follows: consider a weakly increasing
sequence of partitions $P=(\emptyset=\lambda^0,\lambda^1,\dots,\lambda^n)$
where $\lambda^{i-1}$ and $\lambda^i$ differ in size by at most one for
$i\in\{1,2,\dots,n\}$.  (Note that such a sequence corresponds to a partial
Young tableau as described in Section~\ref{sec:RSK-Greene}.)  To $P$ we
associate $jdt(P)=(\emptyset=\mu^0,\mu^1,\dots,\mu^{n-1})$, with the same
property as follows:

\begin{figure}[h]
\begin{equation*}
\def\L#1{\POS[]+(0,5)*{\lambda^{#1}}}
\def\M#1{\POS[]-(0,5)*{\mu^{#1}}}
\begin{xy}
  \xymatrix@!=4pt{%
\e\r\L0&1\dr\L1&2\dr\L2&21\dr\L3&211\dr\L4&211\dr\L5&311\dr\L6&321\dr\L7&3211\dr\L8&3311\d\L9\\
       &\e\r\M0&1 \r\M1&11 \r\M2&111 \r\M3&111 \r\M4&211 \r\M5&221 \r\M6&2211 \r\M7&3211\r\M8&3311\M9}
\end{xy}
\end{equation*}
\caption{jeu de taquin and promotion}
\label{fig:JeuDeTaquin}
\end{figure}

If $\lambda^1=\emptyset$, we set $\mu^i=\lambda^{i+1}$ for $i<n$.  Otherwise,
let $\mu^0=\emptyset$. Suppose that we have already constructed $\mu^{i-1}$ for
some $i<n$.  Then we distinguish three cases: if $\lambda^{i+1}=\lambda^i$, then
we set $\mu^i=\mu^{i-1}$.

If $\nu$ is the only partition that contains $\mu^{i-1}$ and is contained in
$\lambda^{i+1}$, we set $\mu^i=\nu$. Otherwise, there will be exactly one such
partition different from $\lambda^i$, and we set $\mu^i$ equal to this
partition.

The following variation of \jdt is usually called \Dfn{promotion}:
$\overline{jdt}(P)$ is the sequence of partitions obtained by appending the
final partition of $P$ to $jdt(P)$.  An example for this algorithm can be found
in Figure~\ref{fig:JeuDeTaquin}.  Note that promotion is invertible.

\subsection{Jeu de Taquin and Promotion for Semi-Standard Young Tableaux}
\label{sec:jeu-de-taquin-semi}

The description of \jdt is slightly more complicated for semi-standard Young
tableaux, especially in terms of local rules, which are not as well known as in
the standard case.

Given a semi-standard Young tableau $P$, we define $jdt(P)$ as the result of
the following algorithm:
\begin{enumerate}
\item subtract one from all the entries in the tableau.
\item if present, replace the cells containing $0$ by empty cells.  Otherwise
  stop.
\item pick the rightmost empty cell.
\item if there is no cell to the right and no cell above the empty cell, remove
  the empty cell.  Go to Step~3.
\item consider the cells above and to the right of the empty cell.  If there is
  only one such cell, exchange it with the empty cell.  If both cells contain
  the same entry, exchange the upper cell with the empty cell.  Otherwise
  exchange the cell with the smaller entry and the empty cell.  Go to Step~4.
\end{enumerate}
We remark that the conditions in Step~5 are conceived in such a way that the
result is again a semi-standard Young tableau.

We now reproduce the description of \jdt with local rules, as given by Tom
Roby, Frank Sottile, Jeff Stroomer and Julian West in
\cite[Theorem~4.2]{RobySottileStroomerWest2001}.  Given a weakly increasing
sequence of partitions $P=(\emptyset=\lambda^0,\lambda^1,\dots,\lambda^n)$
where $\lambda^{i-1}$ and $\lambda^i$ differ by a horizontal strip, we
associate $jdt(P)=(\emptyset=\mu^0,\mu^1,\dots,\mu^{n-1})$, with the same
property as follows:

If $\lambda^1=\emptyset$, we set $\mu^i=\lambda^{i+1}$ for $i<n$.  Otherwise,
let $\mu^0=\emptyset$. Suppose that we have already constructed $\mu^{i-1}$ for
some $i<n$.  If $\lambda^{i+1}=\lambda^i$, set $\mu^i=\mu^{i-1}$.

Otherwise colour the cells of $\lambda^{i+1}$ that are not in $\lambda^i$
green, and the cells that are in $\lambda^i$ but not in $\mu^{i-1}$ red.  Since
$\lambda^{i+1}$ and $\lambda^i$ differ by a horizontal strip, and so do
$\lambda^i$ and $\mu^{i-1}$, a column contains at most two coloured cells, and
if so, they have different colours.  We now produce a new colouring with the
same property as follows: in each column that contains cells of both colours,
exchange the two cells.  Then, rearrange the coloured cells in each row such
that all red cells are to the right of the green cells.  $\mu^i$ is then the
shape obtained by removing the red cells.

As in the standard case, $\overline{jdt}(P)$ is the sequence of partitions
obtained by appending the final partition of $P$ to $jdt(P)$.  For convenience,
if $P$ is a transposed semi-standard Young tableaux, we define $jdt(P)$ as
$jdt(P^t)^t$ and $\overline{jdt}(P)$ as $\overline{jdt}(P^t)^t$.

\subsection{Jeu de Taquin and Promotion for Fillings}
\label{sec:jeu-de-taquin-1}

The main idea of this article is to define operations analogous to
$\overline{jdt}$ for rectangular fillings:
\begin{dfn}\label{dfn:j}
  Let $\pi$ be a filling of a rectangular polyomino, and let $(P,Q)$ be the
  corresponding pair of sequences of partitions.  Then $\J(\pi)$ is the filling
  corresponding to $(\overline{jdt}(P), Q)$.  An example of this transformation
  applied to a standard filling can be found in Figure~\ref{fig:growth}.
\end{dfn}
Note that, like $\overline{jdt}$, this transformation is invertible.
Furthermore, it is important to keep in mind that in the case of non-standard
fillings the result of $\J$ depends on the variation of the
Robinson-Schensted-Knuth correspondence employed.

In general the transformation $\J$ does not preserve the number of entries of a
given size, if we are using one of the first two methods of
Section~\ref{sec:RS-variations} to standardise the diagram.  For example, using
Burge's method, $\young(11,\hfil1)$ is mapped to $\young(2\hfil,\hfil1)$.
However, there is a notable exception to this failure: $0$-$1$-fillings where
each non-zero entry is the only one in its row or column are mapped to
$0$-$1$-fillings with the same restriction.  Of course, if we use the method
corresponding to RSK' or dual RSK, this is also the case.

The following proposition is a consequence of~\cite[Corollary A.1.2.6]{EC2}, as
pointed out in the proof of~\cite[A.1.2.10]{EC2}:
\begin{prop}\label{prop:jeu-de-taquin}
  Let $\pi$ be a filling of a rectangular polyomino, and let $Q$ be the
  sequence of partitions in the top row of the associated growth diagram.  Let
  $R$ be the sequence of partitions in the top row of the growth diagram
  corresponding to $\pi$ with the first column deleted.  Then $R=jdt(Q)$.
\end{prop}

This proposition might suggest the following alternative description of
$\J(\pi)$: first standardise $\pi$ as in Section~\ref{sec:RS-variations}.
Suppose that the first column of $\pi$ has column-sum $k$, then apply $\J$ as
defined for partial fillings $k$ times.  Finally shrink back the diagram again.

\begin{figure}[h]
  \begin{equation*}
    \begin{xy}<0pt,-42pt>;<7pt,-42pt>:
      (0,0);(12,0)**[|(3)]\dir{-},
      (0,2);(12,2)**\dir{-},
      (0,4);(12,4)**\dir{-},
      (0,6);(12,6)**[|(3)]\dir{-},
      (0,8);(12,8)**\dir{-},
      (0,10);(12,10)**[|(3)]\dir{-},
      (0,12);(12,12)**[|(3)]\dir{-},
      (0,0);(0,12)**[|(3)]\dir{-},
      (2,0);(2,12)**\dir{-},
      (4,0);(4,12)**\dir{-},
      (6,0);(6,12)**[|(3)]\dir{-},
      (8,0);(8,12)**\dir{-},
      (10,0);(10,12)**[|(3)]\dir{-},
      (12,0);(12,12)**[|(3)]\dir{-},
      (1,1),*{\times},(3,3),*{\times},(5,7),*{\times},
      (7,5),*{\times},(9,11),*{\times},(11,9),*{\times},
    \end{xy}
    \quad\mapsto\quad
    \begin{xy}<0pt,-42pt>;<7pt,-42pt>:
      (0,0);(12,0)**[|(3)]\dir{-},
      (0,2);(12,2)**\dir{-},
      (0,4);(12,4)**\dir{-},
      (0,6);(12,6)**[|(3)]\dir{-},
      (0,8);(12,8)**\dir{-},
      (0,10);(12,10)**[|(3)]\dir{-},
      (0,12);(12,12)**[|(3)]\dir{-},
      (0,0);(0,12)**[|(3)]\dir{-},
      (2,0);(2,12)**\dir{-},
      (4,0);(4,12)**[|(3)]\dir{-},
      (6,0);(6,12)**[|(3)]\dir{-},
      (8,0);(8,12)**\dir{-},
      (10,0);(10,12)**\dir{-},
      (12,0);(12,12)**[|(3)]\dir{-},
      (1,1),*{\times},(3,7),*{\times},(5,3),*{\times},
      (7,5),*{\times},(9,11),*{\times},(11,9),*{\times},
    \end{xy}
    \quad\neq\quad
    \begin{xy}<0pt,-42pt>;<7pt,-42pt>:
      (0,0);(12,0)**[|(3)]\dir{-},
      (0,2);(12,2)**\dir{-},
      (0,4);(12,4)**\dir{-},
      (0,6);(12,6)**[|(3)]\dir{-},
      (0,8);(12,8)**\dir{-},
      (0,10);(12,10)**[|(3)]\dir{-},
      (0,12);(12,12)**[|(3)]\dir{-},
      (0,0);(0,12)**[|(3)]\dir{-},
      (2,0);(2,12)**\dir{-},
      (4,0);(4,12)**[|(3)]\dir{-},
      (6,0);(6,12)**[|(3)]\dir{-},
      (8,0);(8,12)**\dir{-},
      (10,0);(10,12)**\dir{-},
      (12,0);(12,12)**[|(3)]\dir{-},
      (1,7),*{\times},(3,11),*{\times},(5,1),*{\times},
      (7,3),*{\times},(9,5),*{\times},(11,9),*{\times},
    \end{xy}
  \end{equation*}
  \caption{standardisation does not commute with promotion}
  \label{fig:j-via-jdt}
\end{figure}
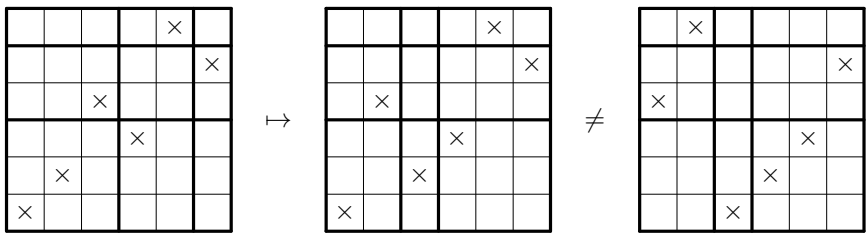

However, it turns out that standardisation does not commute with promotion, and
therefore this idea will not work.  In particular, the non-zero entries in the
last $k$ columns will in general not form a north-east or south-east chain as
they should in the standardisation of a filling.  For example, consider the
filling $\young(\hfil1\hfil,1\hfil1,21\hfil)$, and suppose we are interested in
North-East chains, i.e., we want to use RSK.  Then, the standardised filling is
as on the left of Figure~\ref{fig:j-via-jdt}.  Applying $\J$ three times gives
the filling in the middle of Figure~\ref{fig:j-via-jdt}, which clearly differs
from the standardisation of $\J(\pi)=\young(1\hfil\hfil,1\hfil1,\hfil12)$,
which is depicted on the right.

In Lemma~\ref{lem:commutation} we will see that a slight weakening of such a
commutation property is true.

\begin{dfn}\label{dfn:Knuth-equivalence}
  Two (arbitrary) fillings of a rectangular polyomino are \Dfn{Knuth
    equivalent} if the corresponding partitions labelling the corners along the
  right border are the same.  They are \Dfn{dual Knuth equivalent} if the
  corresponding partitions labelling the corners along the upper border are the
  same.
\end{dfn}

Our proof of Jakob Jonsson's conjecture relies heavily on the following two
propositions, that show in which situation Knuth equivalence is preserved.
\begin{prop}\label{prop:monoid}
  Suppose that the fillings $\young(\alpha)$ and $\young(\alphaP)$ are Knuth
  equivalent, and so are $\young(\beta)$ and $\young(\betaP)$.  Then also
  $\young(\alpha\beta)$ and $\young(\alphaP\betaP)$ are Knuth equivalent.

  Similarly, if $\young(\alpha)$ and $\young(\alphaP)$ are dual Knuth
  equivalent, and so are $\young(\beta)$ and $\young(\betaP)$, then also
  $\young(\beta,\alpha)$ and $\young(\betaP,\alphaP)$ are dual Knuth
  equivalent.
\end{prop}
\begin{proof}
  This is just saying that the plactic monoid is indeed a monoid.
\end{proof}

\begin{prop}\label{prop:chains}
  Consider the filling $\pi^\prime$ defined by columns $i+1,i+2,\dots,i+k$,
  $i\ge 1$, of a filling $\pi$ of a rectangular polyomino.  Then $\pi^\prime$
  is dual Knuth equivalent to the filling defined by columns
  $i,i+1,\dots,i+k-1$ of $\J(\pi)$.  Furthermore, the filling defined by rows
  $i,i+1,\dots,i+k$ of $\pi$ is Knuth equivalent to the filling defined by the
  same rows of $\J(\pi)$.  In particular, row sums are preserved and column
  sums are cyclically shifted by one.
\end{prop}
\begin{proof}
  Let $Q$ be the sequence of partitions labelling the upper border of the
  growth diagram corresponding to $\pi$.  To obtain the sequence of partitions
  along the upper border of $\pi^\prime$, we have to delete the first $i$
  columns of $\pi$ and take the first $k$ partitions labelling the upper border
  of the growth diagram corresponding to the resulting filling.  By
  Proposition~\ref{prop:jeu-de-taquin}, this is equivalent to applying $jdt$
  exactly $i$ times to $Q$ and keeping only the first $k$ partitions.

  On the other hand the filling defined by columns $i, i+1,\dots, i+k-1$ of
  $\J(\pi)$ is obtained by applying $\J$, then deleting the first $i-1$
  columns.  By the definition of $\J$ and Proposition~\ref{prop:jeu-de-taquin},
  the corresponding partitions are again the first $k$ of $jdt$ applied $i$
  times to $Q$.

  To prove the second statement, note first that the sequence of partitions $P$
  along the right border of $\pi$ and $\J(\pi)$ are the same by definition.  To
  obtain the sequence of partitions of the filling defined by rows
  $i,i+1,\dots,i+k$ in $\pi$ or $\J(\pi)$, we can apply $jdt$ exactly $i-1$
  times to $P$ and finally drop all but the first $k+1$ partitions.
\end{proof}

With the aid of the preceding proposition we can also prove the announced
commutation property of standardisation and promotion.  This will also be
needed in the proof of Proposition~\ref{prop:equivalence}.
\begin{lem}\label{lem:commutation}
  Let $\pi$ be an arbitrary filling of a rectangle and $\bar\pi$ the following
  partial standardisation of $\pi$: each column but the first and each row is
  replaced by as many columns and rows as it contains non-zero entries,
  counting multiplicities.  The non-zero entries are then placed in chains
  according to one of the usual standardisation rules.  Then we have
  $std\left(\J\left(\pi\right)\right) =
  std\left(\J\left(\bar\pi\right)\right)$.
\end{lem}
For example, if $\pi$ is $\young(\hfil1\hfil,1\hfil1,21\hfil)$, then, using
RSK-standardisation, $\bar\pi$ is
$\begin{xy}<0pt,-30pt>;<5pt,-30pt>:
  (0,0);(8,0)**[|(3)]\dir{-},
  (0,2);(8,2)**\dir{-},
  (0,4);(8,4)**\dir{-},
  (0,6);(8,6)**[|(3)]\dir{-},
  (0,8);(8,8)**\dir{-},
  (0,10);(8,10)**[|(3)]\dir{-},
  (0,12);(8,12)**[|(3)]\dir{-},
  (0,0);(0,12)**[|(3)]\dir{-},
  (2,0);(2,12)**[|(3)]\dir{-},
  (4,0);(4,12)**\dir{-},
  (6,0);(6,12)**[|(3)]\dir{-},
  (8,0);(8,12)**[|(3)]\dir{-},
  (1,1),*{\times},(1,3),*{\times},(1,7),*{\times},
  (3,5),*{\times},(5,11),*{\times},(7,9),*{\times},
\end{xy}
\,$.
\begin{proof}
  Since $\J$ preserves Knuth equivalence, $std\left(\J\left(\pi\right)\right)$
  and $std\left(\J\left(\bar\pi\right)\right)$ are Knuth equivalent.  We want
  to show that the two fillings are also dual Knuth equivalent.

  By Proposition~\ref{prop:chains} and Proposition~\ref{prop:jeu-de-taquin},
  all but the last $m$ partitions along the upper border of
  $std\left(\J\left(\pi\right)\right)$ and
  $std\left(\J\left(\bar\pi\right)\right)$ are the same, where $m$ is the sum
  of the entries in the first column of $\pi$.

  Furthermore, the but-last partitions of the upper border of $\J(\pi)$ and
  $\J(\bar\pi)$ agree by the definition of $jdt$.  By the definition of
  $\overline{jdt}$, also the final partitions are equal, which implies the
  claim.
\end{proof}

Finally, in Section~\ref{sec:j-globally-well-defined} we will also need the
following proposition, for which, unfortunately, we do not have a short proof.
It asserts that $\J$ preserves locality of modification under a natural
condition.

\begin{prop}\label{prop:equivalence}
  Consider the following two fillings:
  \begin{gather*}
    \young(\alpha\beta\delta)
    \quad\text{and}\quad
    \young(\alpha\gamma\delta)
  \end{gather*}
  and suppose furthermore that $\beta$ and $\gamma$ are Knuth equivalent.
  Then, applying $\J$ to both fillings we obtain
  \begin{gather*}
    \young(\alphaP\betaP\deltaP)
    \quad\text{and}\quad
    \young(\alphaP\gammaP\deltaP)
  \end{gather*}
  where $\alpha^\prime$ has exactly one column less than $\alpha$ and $\delta^\prime$
  has exactly one more column than $\delta$. In this situation, $\beta^\prime$
  and $\gamma^\prime$ are Knuth equivalent.

  Similarly, consider
  \begin{gather*}
    \young(\delta,\beta,\alpha)
    \quad\text{and}\quad
    \young(\delta,\gamma,\alpha)
  \end{gather*}
  and suppose furthermore that $\beta$ and $\gamma$ are dual Knuth equivalent.
  Then, applying $\J$ to both fillings we obtain
  \begin{gather*}
    \young(\deltaP,\betaP,\alphaP)
    \quad\text{and}\quad
    \young(\deltaP,\gammaP,\alphaP)
  \end{gather*}
  where $\alpha^\prime$, $\beta^\prime$ and $\delta^\prime$ have as many rows as
  $\alpha$, $\beta$ and $\delta$ respectively. In this situation, $\beta^\prime$
  and $\gamma^\prime$ are dual Knuth equivalent.
\end{prop}
\begin{proof}
  The proof is given in the appendix.
\end{proof}

\section{Increasing and Decreasing Subsequences in Fillings of Moon
  Polyominoes}
In this section we will apply the transformation $\J$ defined in
Definition~\ref{dfn:j} to moon polyominoes, thus proving a conjecture of Jakob
Jonsson~\cite{Jonsson2005,JonssonWelker2006}.

Various special cases of this theorem were proved by J\"orgen Backelin, Julian
West and Guoce Xin~\cite{BackelinWestXin2004}, Anna de Mier~\cite{DeMier2006},
William Chen, Eva Deng, Rosena Du, Richard Stanley and Catherine
Yan~\cite{ChenDengDuStanleyYan2006}, by Christian
Krattenthaler~\cite{Krattenthaler2006}, and by Jakob Jonsson and Volkmar
Welker~\cite{Jonsson2005,JonssonWelker2006}.  More precisely,
in~\cite{JonssonWelker2006} the special case of stack polyominoes is proved,
using a very different method.  In~\cite{BackelinWestXin2004} and
\cite{ChenDengDuStanleyYan2006} the special case of Ferrers shapes is dealt
with, under various restrictions on the number of non-zero entries allowed,
using completely different methods.  In~\cite{Krattenthaler2006} growth
diagrams were employed to reprove and generalise the results from
\cite{ChenDengDuStanleyYan2006}, and the connection to Jakob Jonsson's
conjecture was noticed.  In the end, the ideas from~\cite{Krattenthaler2006}
led to the proofs given here.  For the precise connection between Christian
Krattenthaler's theorems and the theorem stated below, see
Section~\ref{sec:evacuation}.

\begin{thm}\label{thm:Jakob}
  Let $M$ be a moon polyomino and $\Mat r$ a vector of non-negative integers.
  Let $\Set F_{01}^{ne}(M, l, \Mat r)$ be the set of $0$-$1$-fillings of $M$
  having length of the longest north-east chain equal to $l$ and exactly $\Mat
  r_i$ non-zero entries in row $i$.  Then, for any permutation $\sigma$ of the
  column indices of $M$ such that $\sigma M$ is again a moon polyomino, the
  cardinalities of $\Set F_{01}^{ne}(M, l, \Mat r)$ and $\Set
  F_{01}^{ne}(\sigma M, l, \Mat r)$ coincide.

  In other words, the number of $0$-$1$-fillings of a given moon polyomino with
  given length of longest north-east chain and given number of non-zero entries
  in each row does not depend on the ordering of the columns of the polyomino.

  Furthermore, the number of $0$-$1$-fillings of a given moon polyomino with a
  given number of non-zero entries and given length of longest north-east chain
  depends only on the content of the moon polyomino.
\end{thm}

We prove Theorem~\ref{thm:Jakob} in two steps.  First we show that the
transformation $\J$ from Definition~\ref{dfn:j} can be used to prove a very
general result about cardinalities of unions of chains in arbitrary fillings.
In a second step we show that this implies the theorem above, albeit in a
non-bijective fashion.  Thus, the problem of finding a completely bijective
proof of Theorem~\ref{thm:Jakob} remains open.  However, it appears that this
problem is difficult: Sergi Elizalde~\cite{Elizalde2006} solved the first
non-trivial case, where the polyomino has triangular shape, the length of the
longest north-east chain is two and the number of non-zero entries is maximal.
A stronger bijection for the same situation was recently given by Carlos
Nicol\'as~\cite{Nicolas2009}.

\begin{dfn}
  Let $M$ be a moon polyomino, $\Mat r$ and $\Mat c$ vectors of non-negative
  integers, and $\Lambda$ a mapping that associates to every maximal rectangle
  $R$ in $M$ (which is uniquely determined by its height and its width) a
  partition $\Lambda(R)$.  Then $\Set F^{ne}(M, \Lambda, \Mat r, \Mat c)$ is
  the set of arbitrary fillings of $M$ with
  \begin{itemize}
  \item sum of entries in row $i$ equal to $\Mat r_i$,
  \item sum of entries in column $i$ equal to $\Mat c_i$, and
  \item for every maximal rectangle $R$, the partition
    $\Lambda(R)=(\lambda^{R}_1,\lambda^{R}_2,\dots)$ is such that the maximum
    cardinality of a multiset union of $k$ north-east chains, where any entry
    $e$ in $R$ appears in at most $e$ chains, equals
    $\lambda^{R}_1+\lambda^{R}_2+\dots+\lambda^{R}_k$.
  \end{itemize}
  By the variation of Greene's theorem corresponding to dual RSK',
  Theorem~\ref{thm:Greene-Burge}, the maximum cardinality of a union of $k$
  South-East chains is then given by the sum of the first $k$ parts of the
  transpose of $\Lambda(R)$.

  Let $\Set F^{ne}(M, \Lambda, n)$, be the corresponding set of fillings with
  total sum of entries equal to $n$, but where no attention to row- or column
  sums is paid.

  We also define the analogous sets of fillings for $\Set F^{NE}$, $\Set
  F^{nE}$ and $\Set F^{Ne}$ by replacing \lq north-east\rq\ by one of \lq
  North-East\rq, \lq north-East\rq, or \lq North-east\rq, as appropriate.
  Again, by the appropriate variation of Greene's theorem, we implicitly place
  restrictions on the cardinalities of unions of \lq south-east\rq, \lq
  South-east\rq, and \lq north-East\rq\ chains.
\end{dfn}

\begin{thm}\label{thm:JakobWeak}
  Let $M$ be a moon polyomino and $\sigma$ a permutation of its columns of the
  polyomino, such that $\sigma M$ is again a moon polyomino.  Then the sets of
  maximal rectangles of $M$ and $\sigma M$ coincide and there is an explicit
  bijection between the fillings in $\Set F^{ne}(M, \Lambda, \Mat r, \Mat c)$
  and $\Set F^{ne}(\sigma M, \Lambda, \Mat r, \sigma \Mat c)$, for any value of
  $\Lambda$, $\Mat r$ and $\Mat c$, where $\sigma \Mat c$ is the vector
  obtained by permuting the entries of $\Mat c$ according to $\sigma$.

  Furthermore, the cardinality of $\Set F^{ne}(M, \Lambda, n)$ only depends on
  $\Lambda$, $n$ and the content of $M$, but not $M$ itself.

  Analogous statements hold for $\Set F^{NE}$, $\Set F^{nE}$ and $\Set F^{Ne}$.
\end{thm}
\begin{proof}
  We first show that reordering the columns of the moon polyomino such that the
  result is again a moon polyomino does not change the number of fillings in
  question.  It suffices to show this in the following special case: let $c_1$
  be the $x$-coordinate of any column of the moon polyomino that is contained
  in one of the columns to its right, i.e., a column to the left of the tallest
  column.  Consider the largest rectangle completely contained in the moon
  polyomino that has the same height as the column at $c_1$.  Then we have to
  show that moving the first column of this rectangle to the other end of the
  rectangle, say just after the column with $x$-coordinate $c_2$ does not
  change the number of fillings with the constraints imposed.  For example, we
  could modify a moon polyomino as follows:
  
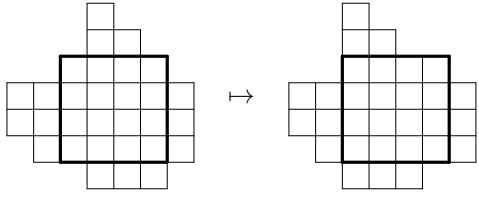
\begin{equation*}
  \begin{xy}<0cm,-35pt>;<10pt,-35pt>:
    (3,0);(6,0)**\dir{-},
    (1,1);(2,1)**\dir{-};(6,1)**[|(3)]\dir{-},(7,1)**\dir{-},
    (0,2);(7,2)**\dir{-},
    (0,3);(7,3)**\dir{-},
    (0,4);(7,4)**\dir{-},
    (2,5);(6,5)**[|(3)]\dir{-},
    (3,6);(5,6)**\dir{-},
    (3,7);(4,7)**\dir{-},
    (0,2);(0,4)**\dir{-},
    (1,1);(1,4)**\dir{-},
    (2,1);(2,5)**[|(3)]\dir{-},
    (3,0);(3,7)**\dir{-},
    (4,0);(4,7)**\dir{-},
    (5,0);(5,6)**\dir{-},
    (6,0);(6,1)**\dir{-};(6,5)**[|(3)]\dir{-},
    (7,1);(7,4)**\dir{-}
  \end{xy}
  \quad\mapsto\quad
  \begin{xy}<0cm,-35pt>;<10pt,-35pt>:
    (2,0);(5,0)**\dir{-},
    (1,1);(2,1)**\dir{-};(6,1)**[|(3)]\dir{-},(7,1)**\dir{-},
    (0,2);(7,2)**\dir{-},
    (0,3);(7,3)**\dir{-},
    (0,4);(7,4)**\dir{-},
    (2,5);(6,5)**[|(3)]\dir{-},
    (2,6);(4,6)**\dir{-},
    (2,7);(3,7)**\dir{-},
    (0,2);(0,4)**\dir{-},
    (1,1);(1,4)**\dir{-},
    (2,0);(2,1)**\dir{-};(2,5)**[|(3)]\dir{-};(2,7)**\dir{-},
    (3,0);(3,7)**\dir{-},
    (4,0);(4,6)**\dir{-},
    (5,0);(5,5)**\dir{-},
    (6,1);(6,5)**[|(3)]\dir{-},
    (7,1);(7,4)**\dir{-}
  \end{xy}
  \end{equation*} 
  In other words, we cyclically shift the columns between $c_1$ and $c_2$.
  
  We now apply the following bijective transformation to the filling of the
  moon polyomino: all the entries outside of the rectangle stay as they are,
  whereas we apply $\J$ to the entries within the rectangle.  We have to show
  that this transformation preserves all the statistics mentioned in the
  statement of the theorem.
  
  By the last statement of Proposition~\ref{prop:chains}, the sum of the
  entries in each row remains the same, and the column sums are cyclically
  shifted, exactly as the columns themselves.

  We first show that there is a one-to-one correspondence between maximal
  rectangles in the two polyominoes.  Consider any maximal rectangle in the
  original polyomino with first column at $d_1$ and last column at $d_2$.
  Then, because the polyomino is intersection-free, we have either $d_1\leq
  c_1< c_2\leq d_2$ or $c_1<d_1\leq d_2\leq c_2$.

  In the first case, the corresponding rectangle, of the same height and width,
  in the new polyomino also consists of the columns between $d_1$ and $d_2$.
  In the second case, because of the cyclic shift, it consists of the columns
  from $d_1-1$ up to $d_2-1$.
  
  It remains to show that the fillings of corresponding maximal rectangles have
  for every $k$ the same maximum cardinality of a union of $k$ north-east
  chains.  (Recall that $\J$ depends on the variation of the
  Robinson-Schensted-Knuth correspondence employed, that is, this determines
  whether we are going to preserve maximal cardinalities of unions of
  north-east, North-East, north-East or North-east-chains.  For brevity, we
  assume that we want to preserve north-east chains here.)

  We argue that, more precisely, if $d_1\leq c_1<c_2\leq d_2$, the fillings of
  the two rectangles are Knuth equivalent, whereas, if $c_1<d_1\leq d_2\leq
  c_2$ they are dual Knuth equivalent.

  To see this, we have to combine Propositions~\ref{prop:chains}
  and~\ref{prop:monoid}: suppose first that $d_1\leq c_1<c_2\leq d_2$, and that
  the original filling of the rectangle is $\young(\alpha\beta\delta)$, such
  that $\beta$ is the restriction of the filling between $c_1$ and $c_2$ to the
  appropriate rows.  The transformation we used leaves entries outside the
  maximal rectangle between $c_1$ and $c_2$ unmodified, therefore, the filling
  of the rectangle in the transformed polyomino is
  $\young(\alpha\betaP\delta)$.  By Proposition~\ref{prop:chains}, $\beta$ and
  $\betaP$ are Knuth equivalent, and therefore, by
  Proposition~\ref{prop:monoid}, $\young(\alpha\beta\delta)$ and
  $\young(\alpha\betaP\delta)$ must be Knuth equivalent, too.

  In case $c_1<d_1\leq d_2\leq c_2$, suppose that the original filling of the
  maximal rectangle between columns $d_1$ and $d_2$ is
  $\young(\delta,\beta,\alpha)$, where $\beta$ is the restriction of the
  maximal rectangle between $c_1$ and $c_2$ to the columns between $d_1$ and
  $d_2$.  Similarly, let $\young(\delta,\betaP,\alpha)$ be the filling of the
  rectangle in the transformed polyomino.  By Proposition~\ref{prop:chains},
  $\beta$ and $\betaP$ are dual Knuth equivalent, and therefore, by
  Proposition~\ref{prop:monoid}, $\young(\delta,\beta,\alpha)$ and
  $\young(\delta,\betaP,\alpha)$ must be dual Knuth equivalent, too.

  To prove the second claim, let $F$ be the Ferrers shape with column heights
  given by the content of $M$.  We exhibit a bijection between fillings in
  $\Set F^{ne}(M, \Lambda, n)$ and $\Set F^{ne}(F, \Lambda, n)$ for any value
  of $\Lambda$ and $n$.

  To do so, we first sort the columns according to their height, using the
  transformation just described, in decreasing order.  This is possible,
  because moon polyominoes are intersection-free.  We then reflect the
  polyomino about the line $x=y$, and obtain a stack polyomino.  Note that
  reflecting the polyomino preserves ne-chains.  Once more we sort the columns
  of the resulting stack polyomino according to height, preserving maximal
  cardinalities of unions of ne-chains.  Finally, we reflect the result again
  about the line $x=y$ and obtain a filling of a Ferrers shape with the same
  content as the original moon polyomino.

  The proof showing that $\Set F^{NE}(M, \Lambda, n)$ and $\Set F^{NE}(F,
  \Lambda, n)$ have the same cardinality is identical to the foregoing,
  however, for the other two cases we have to point out a subtlety: suppose
  that we want to preserve maximal cardinalities of unions of nE-chains.  As
  before, we first sort the columns according to their height transforming the
  filling such that nE-chains are preserved.  Then we reflect the polyomino
  about the line $x=y$.  However, reflection transforms nE-chains into
  Ne-chains.  Therefore, we have to use the variant of $\J$ that preserves
  \emph{Ne} chains to transform the stack-polyomino into a Ferrers shape.  Upon
  reflecting again, we obtain the desired result.
\end{proof}

The proof above does not imply Theorem~\ref{thm:Jakob}: as we observed before,
the transformation $\J$ does not preserve the number of entries of a given
size.  However, a simple inductive argument allows us to reduce it to the
statement about arbitrary fillings of Theorem~\ref{thm:JakobWeak}:
\begin{proof}[Proof of Theorem~\ref{thm:Jakob}]
  Let $\Set F(M, l, \Mat r, \Mat m)$ be the set of arbitrary fillings of $M$
  having length of the longest north-east chain equal to $l$, sum of entries in
  row $i$ equal to $\Mat r_i$ and $\Mat m_i$ non-zero entries in row $i$.  Note
  that $\Set F(M, l, \Mat r, \Mat r)$ is just the set of $0$-$1$-fillings $\Set
  F_{01}^{ne}(M, l, \Mat r)$.

  We will show by induction on the total number $\size{\Mat m}$ of non-zero
  entries in the fillings that $\size{\Set F(M, l, \Mat r, \Mat m)}=\size{\Set
    F(\sigma M, l, \Mat r, \Mat m)}$, for any permutation $\sigma$ such that
  $\sigma M$ is again a moon polyomino.

  When $\size{\Mat m} = 1$, there is only one non-zero entry in the filling,
  and therefore $\Set F(M, l, \Mat r, \Mat m)=\Set F(\sigma M, l, \Mat r, \Mat
  m)$.  Suppose now that $\size{\Mat m} > 1$.  By induction, we have that
  \begin{equation*}
    \sum_{\size{\Mat k} < \size{\Mat m}} \size{\Set F(M, l, \Mat m, \Mat k)}
    =\sum_{\size{\Mat k} < \size{\Mat m}} \size{\Set F(\sigma M, l, \Mat m, \Mat k)}.
  \end{equation*}
  Furthermore, by Theorem~\ref{thm:JakobWeak}, we know that the number of
  arbitrary fillings of $M$ and $\sigma M$ with sum of entries in row $i$ equal
  to $\Mat m_i$ coincide:
  \begin{equation*}
    \sum_{\size{\Mat k} \leq \size{\Mat m}} \size{\Set F(M, l, \Mat m, \Mat k)}
    =\sum_{\size{\Mat k} \leq \size{\Mat m}} \size{\Set F(\sigma M, l, \Mat m, \Mat k)}.
  \end{equation*}
  Therefore we must have equality of $\size{\Set F(M, l, \Mat m, \Mat m)}$ and
  $\size{\Set F(\sigma M, l, \Mat m, \Mat m)}$, i.e., the number of
  $0$-$1$-fillings of the two polyominoes coincide.

  To complete the induction step, it remains to show that 
  \begin{equation*}
    \size{\Set F(M, l, \Mat r, \Mat m)}=\size{\Set F(\sigma M, l, \Mat r, \Mat m)}   
  \end{equation*}
  for arbitrary $\Mat r$.  To do so, let $\Mat k_i$ be a composition of $\Mat
  r_i$ with $\Mat m_i$ parts.  Replacing the ones in a filling in $\Set F(M, l,
  \Mat m, \Mat m)$ -- which is a $0$-$1$-filling -- with the parts of the
  composition, we obtain a filling in $\Set F(M, l, \Mat r, \Mat m)$.
  Moreover, every filling in $\Set F(M, l, \Mat r, \Mat m)$ is obtained in such
  a way.  Since the number of ones in members of $\Set F(M, l, \Mat m, \Mat m)$
  and $\Set F(\sigma M, l, \Mat m, \Mat m)$ coincide in every row, we are done.

  The proof of the second claim, namely that the number of $0$-$1$-fillings of
  a given moon polyomino with a given number of non-zero entries and given
  length of longest north-east chain depends only on the content of the
  polyomino, is quite similar.  Let $\Set F(M, l, n, m)$ be the set of
  arbitrary fillings of $M$ having length of the longest north-east chain equal
  to $l$, sum of entries equal to $n$ and $m$ non-zero entries.  We show by
  induction on $m$ that $\size{\Set F(M, l, n, m)}=\size{\Set F(F, l, n, m)}$,
  where $F$ is the Ferrers shape with column heights given by the content of
  $M$.  The only difference to the foregoing is, that we now have to invoke the
  second statement of Theorem~\ref{thm:JakobWeak}, saying that the the number
  of arbitrary fillings of $M$ only depends on the contents of $M$.
\end{proof}

We conclude this section with counterexamples to a few seemingly natural
generalisations of Theorem~\ref{thm:Jakob}.  First, note that we cannot
restrict the length of the longest NE-chain (or, if preferred: SE-chain)
instead, not even for stack polyominoes.  Indeed, consider the following
filling of a stack polyomino:
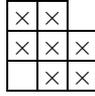
\begin{equation*}
  \young(\times\times,\times\times\times,\hfil\times\times)
\end{equation*}
Its longest NE-chain has length $3$.  However, there is no such filling with
seven non-zero entries of the stack polyomino
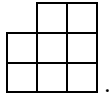
\begin{equation*}
  \Yvcentermath0\Yinterspace0.6ex plus 0.3ex
  \young(:\hfil\hfil,\hfil\hfil\hfil,\hfil\hfil\hfil).
\end{equation*}
In particular, this also implies that one cannot hope to be able to preserve
both the length of the longest ne- and SE-chains.
  
Similarly, we cannot preserve simultaneously the length of the longest ne- and
se-chain.  The following example shows that this is impossible if we
additionally insist on preserving the number of non-zero entries in each row:
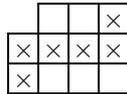
\begin{equation*}
  \young(:\hfil\hfil\times,\times\times\times\times,\times\hfil\hfil\hfil)
\end{equation*}
is a filling with longest north-east chain having length two, and longest
south-east chain having length one.  On the other hand, there is no such
filling of the polyomino
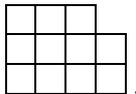
\begin{equation*}
  \Yvcentermath0\Yinterspace0.6ex plus 0.3ex
  \young(\hfil\hfil\hfil:,\hfil\hfil\hfil\hfil,\hfil\hfil\hfil\hfil).
\end{equation*}
Anna de Mier~\cite{DeMier2006} found a counterexample also for the general
case.

Finally, we remark that it is not possible to preserve both row and column sums
analogous to the statement in Theorem~\ref{thm:JakobWeak}.  Consider for
example the filling
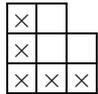
\begin{equation*}
  \young(\times\hfil,\times\hfil\hfil,\times\times\times).
\end{equation*}
There is only one filling with row sums $3,1,1$ and column sums $1,3,1$ of the
permuted polyomino, namely
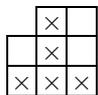
\begin{equation*}
  \young(:\times\hfil,\hfil\times\hfil,\times\times\times),
\end{equation*}
but this has a north-east chain of length two.

\section{A Commutation Property}
\label{sec:j-globally-well-defined}

In this section we would like to prove another beautiful feature of the
transformation $\J$ as applied in the proof of Theorem~\ref{thm:JakobWeak}:
\begin{thm}\label{thm:j-globally-well-defined}
  Applications of $\J$ to different rectangles of a filling of a moon-polyomino
  $M$, such that the first columns of the rectangles are as tall as the
  respective columns of $M$, commute with each other.
\end{thm}
\def\barB{\Bar\beta}
\def\Bprime{\beta^\prime}

  \begin{figure}[h]
    \def\arraystretch{4}
    \begin{tabular}{l@{\extracolsep{8pt}}l}
      \multicolumn{2}{c}{
        \begin{xy}<14pt,0pt>:<0pt,20pt>::
          (0,7)*{\pi},
          (4,0);(8,0)**\dir{-},
          (0,2);(10,2)**\dir{-},
          (0,4);(4,4)**\dir{-},
          (4,4);(8,4)**\dir{-} ?*!U(1.5){U},
          (8,4);(10,4)**\dir{-},
          (4,7);(8,7)**\dir{-},
          (0,2);(0,4)**\dir{-},
          (4,0);(4,7)**\dir{-},
          (8,0);(8,2)**\dir{-},
          (8,2);(8,4)**\dir{-} ?*!R(1.5){T},
          (8,4);(8,7)**\dir{-},
          (10,2);(10,4)**\dir{-},
          (2,3)*{\alpha},(6,3)*{\beta},(9,3)*{\gamma},
          (6,5.5)*{\epsilon},(6,1)*{\delta},
          (1,2);(1,4)**\dir{.},(5,0);(5,7)**\dir{.},
          (0.5,4.5)*{c_1},(4.5,7.5)*{c_2},
        \end{xy}}\\
      \begin{xy}<14pt,0pt>:<0pt,20pt>::
        (0,7)*{\Bar\pi},
        (4,0);(8,0)**\dir{-},
        (1,2);(11,2)**\dir{-},
        (1,4);(4,4)**\dir{-},
        (4,4);(8,4)**\dir{-} ?*!U(1.5){\Bar U},
        (8,4);(11,4)**\dir{-},
        (4,7);(8,7)**\dir{-},
        (1,2);(1,4)**\dir{-},
        (4,0);(4,7)**\dir{-},
        (8,0);(8,2)**\dir{-},
        (8,2);(8,4)**\dir{-} ?*!R(1.5){\Bar T},
        (8,4);(8,7)**\dir{-},
        (11,2);(11,4)**\dir{-},
        (2,3)*{\Bar\alpha},(6,3)*{\Bar\beta},(9,3)*{\Bar\gamma},
        (6,5.5)*{\epsilon},(6,1)*{\delta},
        (10,2);(10,4)**\dir{.},(5,0);(5,7)**\dir{.},
        (10.5,4.5)*{c_1},(4.5,7.5)*{c_2},
      \end{xy}
      &
      \begin{xy}<14pt,0pt>:<0pt,20pt>::
        (0,7)*{\pi^\prime},
        (4,0);(8,0)**\dir{-},
        (0,2);(10,2)**\dir{-},
        (0,4);(4,4)**\dir{-},
        (4,4);(8,4)**\dir{-} ?*!U(1.5){U^\prime},
        (8,4);(10,4)**\dir{-},
        (4,7);(8,7)**\dir{-},
        (0,2);(0,4)**\dir{-},
        (4,0);(4,7)**\dir{-},
        (8,0);(8,2)**\dir{-},
        (8,2);(8,4)**\dir{-} ?*!R(1.25){T^\prime},
        (8,4);(8,7)**\dir{-},
        (10,2);(10,4)**\dir{-},
        (2,3)*{\alpha},(6,3)*{\beta^\prime},(9,3)*{\gamma},
        (6,5.5)*{\epsilon^\prime},(6,1)*{\delta^\prime},
        (1,2);(1,4)**\dir{.},(7,0);(7,7)**\dir{.},
        (0.5,4.5)*{c_1},(7.5,7.5)*{c_2},
      \end{xy}
      \\ 
      \begin{xy}<14pt,0pt>:<0pt,20pt>::
        (0,7)*{\Bar{\Bar\pi}},
        (4,0);(8,0)**\dir{-},
        (1,2);(11,2)**\dir{-},
        (1,4);(4,4)**\dir{-},
        (4,4);(8,4)**\dir{-} ?*!U(1.5){\Bar{\Bar U}},
        (8,4);(11,4)**\dir{-},
        (4,7);(8,7)**\dir{-},
        (1,2);(1,4)**\dir{-},
        (4,0);(4,7)**\dir{-},
        (8,0);(8,2)**\dir{-},
        (8,2);(8,4)**\dir{-} ?*!R(1.5){\Bar{\Bar T}},
        (8,4);(8,7)**\dir{-},
        (11,2);(11,4)**\dir{-},
        (2,3)*{\Bar\alpha},(6,3)*{\Bar{\Bar\beta}},(9,3)*{\Bar\gamma},
        (6,5.5)*{\Bar{\Bar\epsilon}},(6,1)*{\Bar{\Bar\delta}},
        (10,2);(10,4)**\dir{.},(7,0);(7,7)**\dir{.},
        (10.5,4.5)*{c_1},(7.5,7.5)*{c_2},
      \end{xy}
      &
      \begin{xy}<14pt,0pt>:<0pt,20pt>::
        (0,7)*{\pi^{\prime\prime}},
        (4,0);(8,0)**\dir{-},
        (1,2);(11,2)**\dir{-},
        (1,4);(4,4)**\dir{-},
        (4,4);(8,4)**\dir{-} ?*!U(1.5){U^{\prime\prime}},
        (8,4);(11,4)**\dir{-},
        (4,7);(8,7)**\dir{-},
        (1,2);(1,4)**\dir{-},
        (4,0);(4,7)**\dir{-},
        (8,0);(8,2)**\dir{-},
        (8,2);(8,4)**\dir{-} ?*!R(1){T^{\prime\prime}},
        (8,4);(8,7)**\dir{-},
        (11,2);(11,4)**\dir{-},
        (2,3)*{\alpha^{\prime\prime}},(6,3)*{\beta^{\prime\prime}},
        (9,3)*{\gamma^{\prime\prime}},
        (6,5.5)*{\epsilon^\prime},(6,1)*{\delta^\prime},
        (10,2);(10,4)**\dir{.},(7,0);(7,7)**\dir{.},
        (10.5,4.5)*{c_1},(7.5,7.5)*{c_2},
      \end{xy}
    \end{tabular}
    \caption{applying $\J$ successively to $\Gamma$ and $\Delta$}
    \label{fig:rectangles}
  \end{figure}
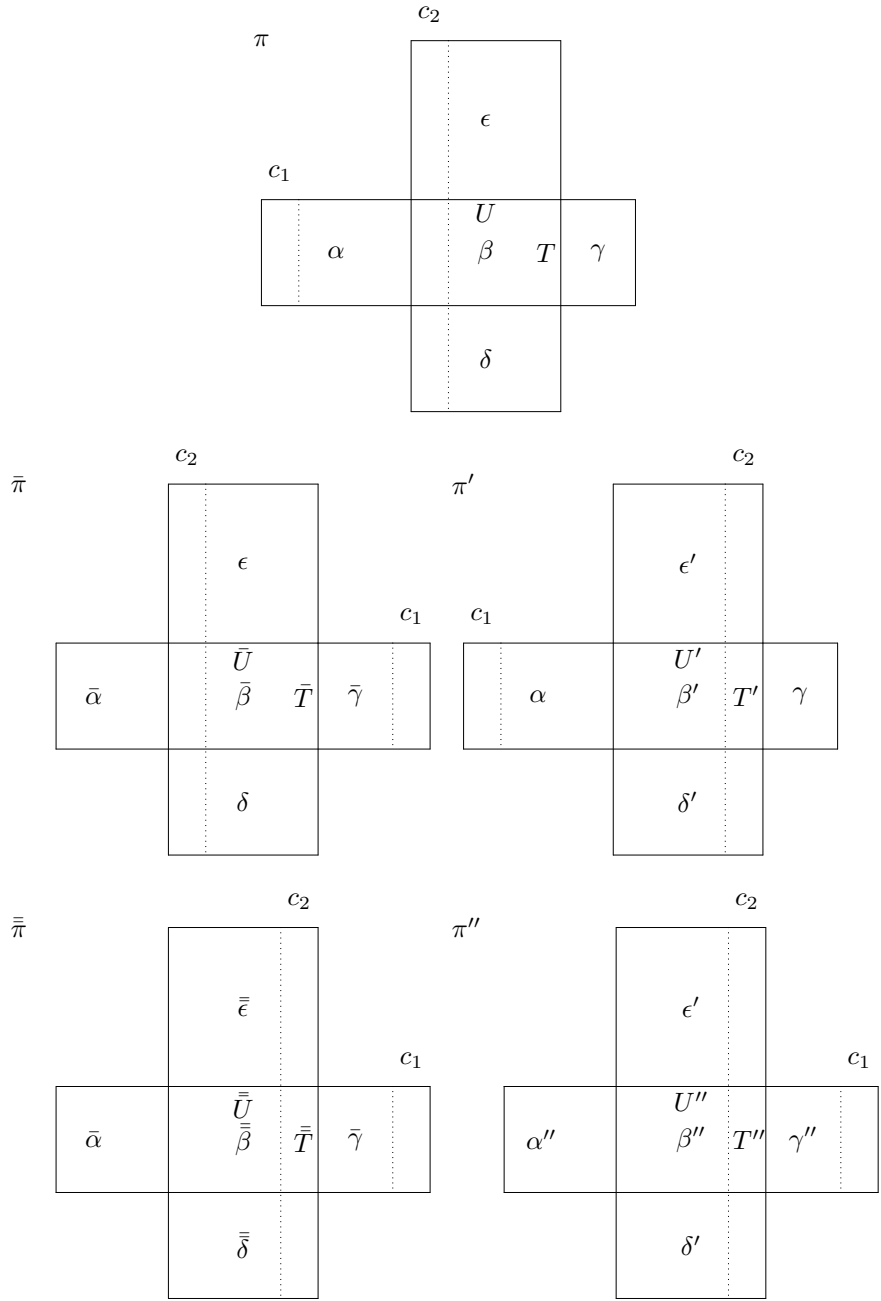
\begin{proof}\Yboxdim13pt
  Suppose that $c_1$ and $c_2$ are respectively the first columns of two
  maximal rectangles $\Gamma$ and $\Delta$ and that we want to move $c_1$ to
  the end of the first rectangle, and $c_2$ towards the end of the second
  rectangle.  Since all entries outside of these two rectangles stay unchanged,
  we can assume that the moon polyomino consists only of the union of $\Gamma$
  and $\Delta$, as schematically depicted at the top of
  Figure~\ref{fig:rectangles}.  There, we subdivided $\Gamma$ into three
  smaller rectangles $\alpha$, $\beta$ and $\gamma$ and $\Delta$ into $\delta$,
  $\beta$ and $\epsilon$.

  In the following we will have to consider several growth diagrams
  simultaneously.
  In particular, we consider the growth diagram corresponding to $\beta$ and
  will label its top and right border with sequences of partitions $U$ and $T$.
  
  On the left hand side of Figure~\ref{fig:rectangles} we see what happens to
  the original filling when we apply $\J$ first to the rectangle
  $\Gamma=\young(\alpha\beta\gamma)$, to obtain $\Bar\pi$, and then to
  $\young(\epsilon,\barB,\delta)$ to obtain $\Bar{\Bar\pi}$. On the right hand
  side the result $\pi^\prime$ of applying $\J$ first to the rectangle
  $\Delta=\young(\epsilon,\beta,\delta)$ and then the result
  $\pi^{\prime\prime}$ of applying $\J$ to $\young(\alpha\Bprime\gamma)$ is
  shown. We have to prove that the fillings $\Bar{\Bar\pi}$ and
  $\pi^{\prime\prime}$ at the bottom of Figure~\ref{fig:rectangles} are the
  same.

  We first observe that $\young(\beta)$ and $\young(\Bprime)$ are Knuth
  equivalent.  This follows by applying the second part of
  Proposition~\ref{prop:chains} to the rows corresponding to $\beta$ in
  $\young(\epsilon,\beta,\delta)$.  Thus we can apply the first part of
  Proposition~\ref{prop:equivalence} to $\young(\alpha\beta\gamma)$ and
  $\young(\alpha\betaP\gamma)$ and obtain that
  $\Bar\alpha=\alpha^{\prime\prime}$, $\Bar\gamma=\gamma^{\prime\prime}$ and
  $\Bar T=T^{\prime\prime}$.  Again by the second part of
  Proposition~\ref{prop:chains}, applied to the rows corresponding to
  $\Bar\beta$ in $\young(\epsilon,\barB,\delta)$ we have that $\Bar T =
  \Bar{\Bar T}$.

  Very similarly, applying the first part of Proposition~\ref{prop:chains}, we
  observe that $\young(\beta)$ and $\young(\barB)$ are dual Knuth equivalent.
  Thus, the second part of Proposition~\ref{prop:equivalence} applied to
  $\young(\epsilon,\beta,\delta)$ and $\young(\epsilon,\barB,\delta)$, shows
  that $\Bar{\Bar\delta}=\delta^\prime$, $\Bar{\Bar\epsilon}=\epsilon^\prime$
  and $\Bar{\Bar U}= U^\prime$.  By the first part of
  Proposition~\ref{prop:chains} we have $U^\prime = U^{\prime\prime}$.

  Finally, $\Bar{\Bar\beta}=\beta^{\prime\prime}$ follows from $\Bar{\Bar
    T}=T^{\prime\prime}$ and $\Bar{\Bar U}=U^{\prime\prime}$.
\end{proof}

\section{Evacuation and Promotion for Stack Polyominoes}
\label{sec:evacuation}
In this section we relate our bijection to evacuation, and thereby to the
construction employed by Christian Krattenthaler~\cite{Krattenthaler2006} to
prove Theorem~\ref{thm:Jakob} and Theorem~\ref{thm:JakobWeak} for the special
case of Ferrers diagrams.  Let us first recall the definition of evacuation.
\begin{dfn}
  Given a weakly increasing sequence of partitions
  $Q=(\emptyset=\lambda^0,\lambda^1,\dots,\lambda^n)$, we construct the
  \Dfn{evacuated} sequence of partitions
  $ev(Q)=(\emptyset=\mu^0,\mu^1,\dots,\mu^n)$ as follows: We set
  $\mu^n=\lambda^n$, and then $\mu^{n-i}$ equal to the last partition of
  $jdt\left(\cdots jdt(Q)\right)$, where we apply $jdt$ $i$ times.
\end{dfn}

Christian Krattenthaler \cite{Krattenthaler2006} used the following bijection
on Ferrers shapes:
\begin{dfn}\label{dfn:e}
  Let $\pi$ be a filling of a Ferrers shape and $\Delta$ the associated growth
  diagram.  Let $e(\Delta)$ be the growth diagram obtained from $\Delta$ by
  transposing all the partitions along the top and right border and applying
  the backward rules B1 to B4 to obtain the remaining partitions and the
  entries of the squares.  Let $e(\pi)$ be the filling associated to
  $e(\Delta)$.
\end{dfn}
We will show that this bijection is to evacuation what our transformation $\J$
is to promotion.  To this end we extend the notion of growth diagrams
introduced in Section~\ref{sec:growth} to stack polyominoes.

\subsection{Growth Diagrams for Stack Polyominoes}
Throughout this section, we fix a variant of standardisation.

\begin{dfn}\label{dfn:stack-growth}
  Given a filling $\pi$ of a stack polyomino, we label the top-right corners of
  the polyomino with \emph{tuples} of partitions to obtain a \Dfn{growth
    diagram} as follows:
  \begin{itemize}
  \item if there is no corner directly above the corner $c$, or the rows just
    above and just below it are left-justified, the label consists of a single
    partition.  This partition is computed according to the appropriate variant
    of Greene's theorem applied to the filling restricted to the largest
    rectangle below and to the left of $c$.
  \item Otherwise, suppose that the row above $c$ is indented to the right by
    $m$ columns with respect to the row just below $c$: the row just below $c$
    begins in column $i$, and the row just below $c$ in column $i-m$. Then the
    label consists of a tuple of $m+1$ partitions: the $j$\textsuperscript{th}
    partition is computed applying Greene's theorem to the filling restricted
    to the largest rectangle below and to the left of $c$, beginning in column
    $i-j+1$.
  \end{itemize}
  We say that the sequence of partitions, read beginning at the left-most
  corner of the top-row, down to the corner at the bottom-right
  \Dfn{corresponds} to $\pi$.
\end{dfn}

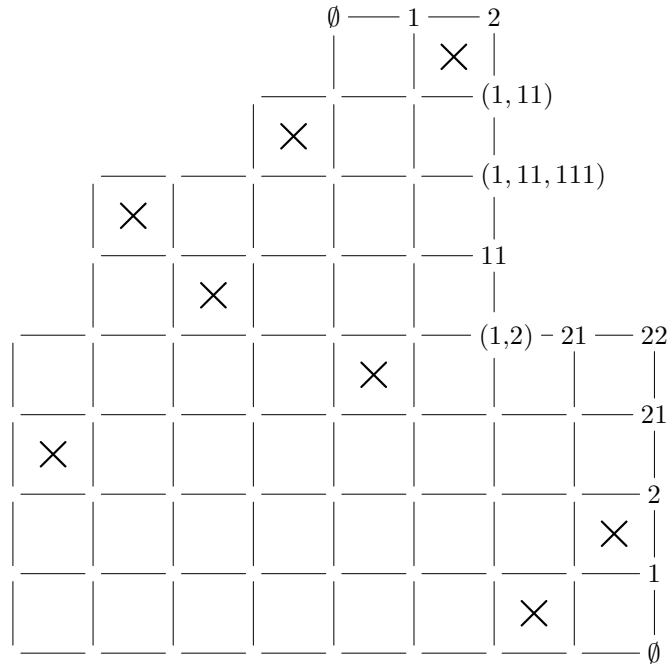
\begin{figure}
  \begin{equation*}
  \begin{xy}
  \xymatrix@!=6pt{%
          &       &       &       &\e\dr&1\dr\y&2  \d&       &      \\
          &       &       &{}\dr\y&{}\dr&{}\dr &**[r]{(1,11)}\d
                                                     &       &      \\
          &{}\dr\y&{}\dr  &{}\dr  &{}\dr&{}\dr &**[r]{(1,11,111)}\d
                                                     &       &\\
          &{}\dr  &{}\dr\y&{}\dr  &{}\dr&{} \dr&11 \d&       &\\
     {}\dr&{}\dr  &{}\dr  &{}\dr  &{}\dr\y&{}\dr&**[r]\raisebox{-9pt}{(1,2)}\dr&21\dr&22\d\\
   {}\dr\y&{}\dr  &{}\dr  &{}\dr  &{}\dr&{} \dr&{}\dr&{} \dr&21  \d\\
   {}\dr  &{}\dr  &{}\dr  &{}\dr  &{}\dr&{} \dr&{}\dr&{}\dr\y&2  \d\\
   {}\dr  &{}\dr  &{}\dr  &{}\dr  &{}\dr&{} \dr&{}\dr\y&{}\dr&1  \d\\
   {}\r   &{}\r   &{}\r   &{}\r   &{}\r &{} \r &{}\r &{}  \r &\e}
  \end{xy}
  \end{equation*}
  \caption{a growth diagram for a stack polyomino}
  \label{fig:stackgrowth}
\end{figure}
An example of such a generalised growth diagram is given in
Figure~\ref{fig:stackgrowth}.  It is clear that for Ferrers shapes the
construction above coincides with the obvious extension of growth diagrams as
presented in Section~\ref{sec:growth} and introduced by Sergey Fomin and Tom
Roby~\cite{Fomin1986,Roby1991}.

Similar to the case of growth diagrams for rectangular shapes we can describe
precisely which sequences of partitions labelling the upper-right border
actually occur.  For partial fillings we have:
\begin{prop}
  Let $S$ be a stack polyomino.  Then there is a bijection between partial
  fillings of $S$ and labellings of the top-right corners of $S$ with tuples of
  partitions satisfying the following conditions:
  \begin{itemize}
  \item the left-most corner in the top row and the corner on the bottom right
    are labelled with the empty partition.
  \item writing tuples of partitions from left to right, if $\mu$ is the
    partition just to the left of $\lambda$, then $\mu$ is obtained from
    $\lambda$ by deleting at most one from some part.
  \item if $\mu$ is the partition just below $\lambda$, then $\mu$ is obtained
    from $\lambda$ by deleting at most one from some part.
  \end{itemize}
\end{prop}
\begin{rmk}
  To obtain correspondences for arbitrary fillings, we only have to replace \lq
  deleting at most one from some part\rq\ by \lq deleting a horizontal
  strip\rq, or \lq deleting a vertical strip\rq, according to the chosen
  variant of standardisation.
\end{rmk}
\begin{proof}
  For left justified rows we can reconstruct the filling from top to bottom
  given the partitions along the top and the right corners using the local
  rules B1 to B4.  Thus, the only subtlety is how to deal with successive rows
  that are not left justified.

  Consider the following situation:
  \begin{figure}[h]
    \begin{equation*}
      \begin{xy}<10pt,0pt>:
        (0,3);(7,3)**\dir{-},
        (0,3.5);(7,3.5)**\dir{}?;(0.1,3.5).(6.9,3.5)**\frm{^\}};*!D(3){m},
        (7,5);(16,5)**\dir{-},
        (7,3);(7,5)**\dir{-},
        (0,0);(0,3)**\dir{-},
        (16,0);(16,5)**\dir{-},
        (8.5,1.5)*{\pi}, 
        (8,3)*{\mu_1},(9.25,3)*{\mu_2},(11.5,3)*{\dots},(15,3)*{\mu_l},
        (16.75,3.2)*[r]{(\rho_0,\rho_1,\dots,\rho_m)},
        (16.75,0)*[r]{(\nu_0,\nu_1,\dots)},
        (1,0)*{\lambda_1},(2.25,0)*{\lambda_2},(8,0)*{\dots},(14.4,0)*{\lambda_{m+l}},
        (0,0);(0.25,0)**\dir{-},
        (2.9,0);(7,0)**\dir{-},
        (9,0);(13,0)**\dir{-},
        (15.75,0);(16,0)**\dir{-},
      \end{xy}
    \end{equation*}
    \caption{reconstructing the partitions labelling the next row}
    \label{fig:lowerpartitions}
  \end{figure}
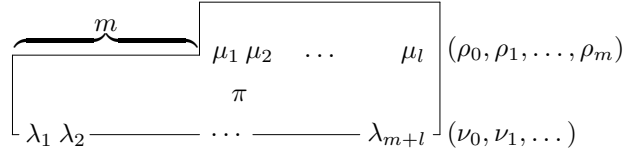
  We are given a tuple $(\rho_0,\rho_1,\dots,\rho_m)$ of partitions labelling
  the top-right corner of the row, another tuple $(\nu_0,\nu_1,\dots)$
  labelling its bottom-right corner, and partitions $\e=\mu_0,\mu_1,
  \dots,\mu_l,\mu_{l+1}=\rho_0$ labelling the corners between the two rows.  We
  would like to determine the filling $\pi$, and the sequence of partitions
  $\e=\lambda_0,\lambda_1, \dots,\lambda_{m+l}$.

  We want to reconstruct the sequence of partitions $Q=(\e=\mu'_0,\mu'_1,
  \dots,\mu'_{m+l})$ such that the last partition of $Q$ is $\rho_m$, the last
  partition of $jdt(Q)$ is $\rho_{m-1}$ and so on.  Finally, applying $jdt$ $m$
  times to $Q$ should yield a sequence of partition ending with $\rho_0$.  By
  Proposition~\ref{prop:jeu-de-taquin}, this is exactly how the partitions
  labelling the corners of the stack polyomino are computed in
  Definition~\ref{dfn:stack-growth}.  Given $Q$ and $\nu_0$ we can then use the
  local rules B1 to B4 to compute the filling $\pi$ and
  $\e=\lambda_0,\lambda_1, \dots,\lambda_{m+l}$.

  By considering the growth diagram description of \jdt given in
  Section~\ref{sec:jeu-de-taquin}, we see that the preimage of $jdt$ can indeed
  be computed, if only one is given the last partition of the preimage.  In the
  first step we compute the preimage of
  $(\e=\mu_0,\mu_1,\dots,\mu_{l+1}=\rho_0)$, assuming as its last partition
  $\rho_1$.  We then compute the preimage of this sequence, assuming as its
  last partition $\rho_2$, and so on.
\end{proof}

\subsection{Evacuation and Promotion}
In this section, we will analyse the effect of reversing the order of the
columns of a filling on the sequence of partitions it corresponds to.  For
rectangular shapes, this is well known:
\begin{prop}[Corollary~A 1.2.11 of \cite{EC2}]
  \label{prop:reverse}
  If $\pi$ corresponds to $(P,Q)$ then the filling obtained by reversing the
  order of the columns of $\pi$ corresponds to $(P^t,ev(Q)^t)$.
\end{prop}
To relate our construction to Christian Krattenthaler's bijection, we will use
the following transformation:
\begin{dfn}
  Let $\pi$ be a filling of a Ferrers shape $F$.  Let $ev^t(\pi)$ be the
  filling of the reversed polyomino obtained by applying $\J^{-1}$ to the
  largest rectangle contained in $F$ spanning all columns, then $\J^{-1}$ to
  the largest rectangle spanning all columns of the result but the first, and
  so on.
\end{dfn}
We can now state and prove the main theorem of this section:
\begin{thm}
  Let $\pi$ be a filling of a Ferrers shape.  Suppose that the sequence of
  partitions labelling its top-left corners, reading from the left of the top
  row down to the bottom-right is
  $(\e=\lambda_0,\lambda_1,\dots,\lambda_n=\e)$, and that its top-row consists
  of $k$ cells.  Let $\pi^r$ be the filling (of a stack-polyomino) obtained by
  reflecting $\pi$ about a vertical line.  Then we have:
  \begin{enumerate}
  \item Let $(\e=\mu_0,\mu_1,\dots,\mu_n=\e)$ be the sequence of partitions
    labelling $\pi^r$.  Then
    $(\mu_0^t,\mu_1^t,\dots,\mu_k^t)=ev(\lambda_0,\lambda_1,\dots,\lambda_k)$,
    and $\mu_i^t=\lambda_i$ for $i\geq k$.
  \item Let $(\e=\nu_0,\nu_1,\dots,\nu_n=\e)$ be the sequence of partitions
    labelling $ev^t(\pi)$.  Then
    $(\nu_0,\nu_1,\dots,\nu_k)=ev(\lambda_0,\lambda_1,\dots,\lambda_k)$, and
    $\nu_i=\lambda_i$ for $i\geq k$.
  \end{enumerate}
\end{thm}
\begin{proof}
  For rectangular shapes, the first statement is precisely
  Proposition~\ref{prop:reverse}.  To prove the second statement for
  rectangular shapes, let the sequence of partitions labelling the right
  corners of $\pi$ be $P$ and let $Q$ be the sequence of partitions labelling
  the top corners.  Then, by Proposition~\ref{prop:jeu-de-taquin}, applying
  $ev^t$ to $\pi$ amounts to applying evacuation on $Q$ and leaving $P$
  unchanged.

  To see the first statement for stack polyominoes, consider corresponding
  partitions in $\pi$ and its reverse, as in Figure~\ref{fig:reverse}.  Again,
  the fillings within the dotted rectangles are just the reversed of each
  other, and Proposition~\ref{prop:reverse} implies $\mu_i=\lambda_i^t$.
  \begin{figure}[h]
    \begin{equation*}
      \begin{xy}<0cm,-35pt>;<10pt,-35pt>:
        (0,0);(7,0)**\dir{-},
        (4,4);(7,4)**\dir{-},
        (0,6);(4,6)**\dir{-},
        (0,0);(0,6)**\dir{-},
        (4,4);(4,6)**\dir{-},
        (7,0);(7,4)**\dir{-},
        (0.1,0.1);(4.9,0.1)**\dir{.},
        (0.1,0.1);(0.1,3.9)**\dir{.},
        (4.9,0.1);(4.9,3.9)**\dir{.},
        (0.1,3.9);(4.9,3.9)**\dir{.},
        (5.4,4.5)*{\lambda_i}
      \end{xy}
      \quad\leftrightarrow\quad
      \begin{xy}<0cm,-35pt>;<10pt,-35pt>:
        (0,0);(7,0)**\dir{-},
        (3,6);(7,6)**\dir{-},
        (0,0);(0,4)**\dir{-},
        (0,4);(3,4)**\dir{-},
        (3,4);(3,6)**\dir{-},
        (7,0);(7,6)**\dir{-},
        (2.1,0.1);(6.9,0.1)**\dir{.},
        (2.1,0.1);(2.1,3.9)**\dir{.},
        (6.9,0.1);(6.9,3.9)**\dir{.},
        (2.1,3.9);(6.9,3.9)**\dir{.},
        (7.6,4.4)*{\mu_i}
      \end{xy}
    \end{equation*}
    \caption{corresponding labels of $\pi$ and its reverse}
    \label{fig:reverse}
  \end{figure}
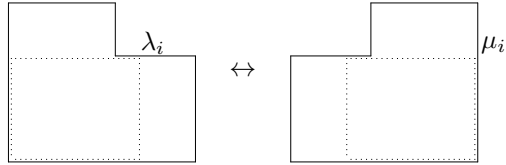

  For the proof of the second statement, suppose that we have already applied
  $\J^{-1}$ several times, and we are now going to apply it to the next block
  of rectangles of common height, as on the left of Figure~\ref{fig:jstar}.

  \begin{figure}[h]
    \begin{equation*}
      \begin{xy}<0cm,-35pt>;<10pt,-35pt>:
        (0,0);(12,0)**\dir{-},
        (0,0);(0,1)**\dir{-},  
        (0,1);(2,1)**\dir{-},  
        (2,1);(2,2)**\dir{-},
        (2,2);(3,2)**\dir{-},
        (3,2);(3,6)**\dir{-},
        (3,6);(6,6)**\dir{-},
        (6,6);(6,5)**\dir{-},
        (6,5);(8,5)**\dir{-},
        (8,5);(8,4)**\dir{-},
        (8,4);(12,4)**\dir{-},
        (12,4);(12,0)**\dir{-},
        (3.1,0.1);(11.9,0.1)**\dir{.},
        (3.1,0.1);(3.1,3.9)**\dir{.},
        (3.1,3.9);(11.9,3.9)**\dir{.},
        (11.9,3.9);(11.9,0.1)**\dir{.},
        (8.6,4.5)*{\lambda_0},
        (9.6,4.5)*{\lambda_1},
        (10.8,4.5)*{\dots},
        (12.6,4.5)*{\lambda_m},
        (7,2)*{\pi'}
      \end{xy}
      \quad\mapsto\quad
      \begin{xy}<0cm,-35pt>;<10pt,-35pt>:
        (0,0);(12,0)**\dir{-},
        (0,0);(0,1)**\dir{-},  
        (0,1);(2,1)**\dir{-},  
        (2,1);(2,2)**\dir{-},
        (2,2);(3,2)**\dir{-},
        (3,2);(3,4)**\dir{-},
        (3,4);(6,4)**\dir{-},
        (6,4);(6,6)**\dir{-},
        (6,6);(9,6)**\dir{-},
        (9,6);(9,5)**\dir{-},
        (9,5);(12,5)**\dir{-},
        (12,5);(12,0)**\dir{-},
        (3.1,0.1);(11.9,0.1)**\dir{.},
        (3.1,0.1);(3.1,3.9)**\dir{.},
        (3.1,3.9);(11.9,3.9)**\dir{.},
        (11.9,3.9);(11.9,0.1)**\dir{.},
        (15.6,4.3)*{(\nu_0,\nu_1,\dots,\nu_m)}
      \end{xy}
    \end{equation*}
    \caption{applying $\J^{-1}$ to a block of rectangles}
    \label{fig:jstar}
  \end{figure}
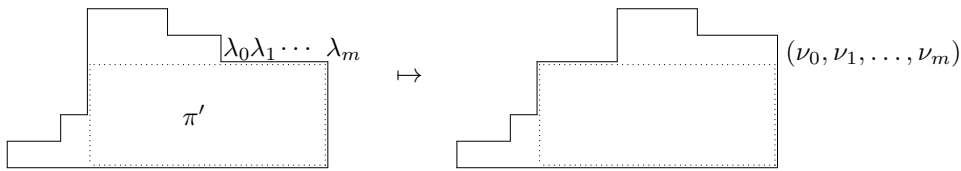

  Since $\J^{-1}$ preserves Knuth-equivalence, and we do not modify the filling
  to the left of the dotted rectangle anymore, Proposition~\ref{prop:monoid}
  implies that the labels below the top row of the dotted rectangle will not be
  altered.  Furthermore, because of Proposition~\ref{prop:chains} and since we
  do not modify the filling above the dotted rectangle when applying $\J^{-1}$
  the next $m$ times, the labels above the top row of the dotted rectangle will
  also stay the same.

  Thus, it remains to show that, with the notation indicated in the picture,
  $\lambda_i=\nu_i$.  To this end, consider the filling $\pi'$ on the left of
  Figure~\ref{fig:jstar} within the dotted rectangle.  We need to compute the
  effect of applying $\J^{-1}$ a total of $m$ times, first to the full
  rectangle, then to all but the first column of the result and so on.

  The first application of $\J^{-1}$ preserves the shape $\lambda_m$ of the
  filling.  In the following, the first column of the resulting filling is not
  modified anymore.  Since $\J^{-1}$ preserves Knuth-equivalence by definition,
  Proposition~\ref{prop:monoid} implies that the shape of the filling of the
  full rectangle is not modified at all during the remaining $m-1$ applications
  of $\J^{-1}$.  Therefore, $\nu_m=\lambda_m$.

  Furthermore, by Proposition~\ref{prop:chains}, the filling obtained from
  $\pi'$ by deleting the last column is dual Knuth equivalent to the filling
  $\pi''$ obtained from $\J^{-1}(\pi')$ by deleting its first column.  In
  particular, the shape of $\pi''$ is $\lambda_{m-1}$, and as before we deduce
  that $\nu_{m-1}=\lambda_{m-1}$.  The claim follows by induction.
\end{proof}

\appendix
\section{Locality of Promotion}
In this appendix we prove Proposition~\ref{prop:equivalence}.  We will first
consider only the case of standard fillings, and then deduce the general case
using standardisation, using Lemma~\ref{lem:commutation}.

\begin{dfn}
  Two partial fillings differ by a \Dfn{Knuth relation of the first kind}, if
  all but three columns coincide, and the remaining three columns are related
  (schematically) by transforming
  \begin{equation*}
    \begin{xy}<0pt,-6pt>;<7pt,-6pt>:
      (0,2)*{\times}="a",(1,3)*{\times};(2,0)*{\times}."a"**\frm{-}
    \end{xy}
    \quad\text{into}\quad
    \begin{xy}<0pt,-6pt>;<7pt,-6pt>:
      (0,2)*{\times}="a",(1,0)*{\times};(2,3)*{\times}."a"**\frm{-}
    \end{xy}\,,
  \end{equation*}
  or vice versa.  They differ by a \Dfn{Knuth relation of the second kind}, if
  the remaining three columns are related by transforming
  \begin{equation*}
    \begin{xy}<0pt,-6pt>;<7pt,-6pt>:
      (0,0)*{\times}="a",(1,3)*{\times};(2,1)*{\times}."a"**\frm{-}
    \end{xy}
    \quad\text{into}\quad
    \begin{xy}<0pt,-6pt>;<7pt,-6pt>:
      (0,3)*{\times}="a",(1,0)*{\times};(2,1)*{\times}."a"**\frm{-}
    \end{xy}\,,
  \end{equation*}
  or vice versa.  The indicated transformations are called \Dfn{Knuth
    transformations}.
\end{dfn}

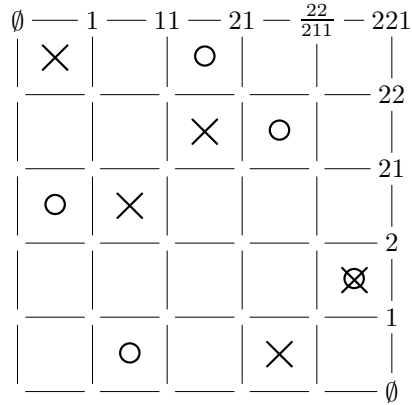
\begin{figure}[h]
\begin{equation*}
\begin{xy}
  \xymatrix@!=4pt{%
  \e\dr\y&1 \dr  &11\dr\c&21\dr  &\frac{22}{211}\dr    &221\d\\
  {}\dr  &{}\dr  &{}\dr\y&{}\dr\c&{}\dr                &22 \d\\
  {}\dr\c&{}\dr\y&{}\dr  &{}\dr  &{}\dr                &21 \d\\
  {}\dr  &{}\dr  &{}\dr  &{}\dr  &{}\dr            \c\y&2  \d\\
  {}\dr  &{}\dr\c&{}\dr  &{}\dr\y&{}\dr                &1  \d\\
  {}\r   &{}\r   &{}\r   &{}\r   &{}\r                 &\e}
\end{xy}
\end{equation*}
\caption{Two very different Knuth equivalent fillings, whose growth diagrams
  differ along the border only in a single partition.  The filling with circles
  is not shape equivalent to a triangle in the last three columns.}
\label{fig:counter-example}
\end{figure}

It is well known that two partial fillings are Knuth equivalent (see
Definition~\ref{dfn:Knuth-equivalence}) if and only if they can be transformed
one into the other using a sequence of Knuth transformations.  The following
proposition gives also some information on the sequence of partitions along the
upper border of the growth diagram.  For its precise statement, we need yet
another definition:
\begin{dfn}
  A partial filling of a rectangle with three columns is \Dfn{shape equivalent
    to a triangle}, if the partition labelling the top right corner of the
  corresponding growth diagram equals $21$.

  In terms of tableaux, we say that three consecutive entries $k-1$, $k$ and
  $k+1$ of a partial tableau are \Dfn{shape equivalent to a triangle} if
  applying $jdt$ $k-2$ times to the tableau results in a tableau where the
  entries $1$, $2$ and $3$ are arranged as $\young(2,13)$ or $\young(3,12)$.
\end{dfn}
\begin{prop}\label{prop:difference}
  Two partial fillings differ by a single Knuth relation (of either kind),
  namely in columns $k-1$, $k$ and $k+1$ if and only if the fillings are Knuth
  equivalent, the sequences of partitions along the upper border of their
  growth diagrams differ in exactly one partition, namely either between
  columns $k-1$ and $k$ or $k$ and $k+1$, and each of the two fillings
  restricted to these three columns is shape equivalent to a triangle.
\end{prop}
Note that there are Knuth equivalent growth diagrams that differ in exactly one
partition along the upper border, but cannot be transformed one into the other
by a single Knuth transformation.  A small example is given in
Figure~\ref{fig:counter-example}, where we superimposed the two growth
diagrams, using crosses for one and circles for the other.  We write the
partitions belonging to the diagram with circles just above the partitions
belonging to the diagram with crosses, wherever the partitions labelling a
corner differ. We will keep this notation throughout the appendix.

\begin{proof}
  \textcircled{\footnotesize$\Rightarrow$} Suppose that the two growth diagrams
  differ by a single Knuth relation of the first kind, as shown in
  Figure~\ref{fig:Knuth-columns}.a. In this case, all partitions along the
  upper border up to and including the partition just after column $k-1$ must
  be equal, since, up to this point the fillings are identical. Furthermore,
  all partitions along the top corner starting with the partition just after
  column $k+1$ must be equal, since the fillings up to that column are Knuth
  equivalent.  We conclude that, in this case, exactly the partitions of the
  two growth diagrams on the top border just after column $k$ are different.

  It remains to check the case where the two growth diagrams differ by a single
  Knuth relation of the second kind. 

  To ease the description, we introduce the following notation: given a
  partition $\lambda$, let $\lambda+\epsilon_k$ be the result of adding one to
  the $k$\textsuperscript{th} part of $\lambda$. Furthermore, we write
  $\lambda+1$ as shorthand for $\lambda+\epsilon_1$. We then have a situation
  as shown schematically in Figure~\ref{fig:Knuth-columns}.b.

  \begin{figure}[h]
  \begin{tabular}{cc}
    \begin{xy}
  \xymatrix@!=10pt{%
{}\dd        &{}\dd                      &{}\dd                 &{}\dd\\
{}\dr        &{}\dr \y                   &{}\dr          \c     &{}\d\\
{}\ddr       &{}\ddr                     &{}\ddr                &{}\dd\\
{}\dr    \y\c&{}\dr                      &{}\dr                 &{}\d\\
{}\ddr       &{}\ddr                     &{}\ddr                &{}\dd\\
{}\dr        &{}\dr                    \c&{}\dr          \y     &{}\d\\
{}\ddr       &{}\ddr                     &{}\ddr                &{}\dd\\
{}\put{$k-1$}&{}\put{$k$}                &{}\put{$k+1$}         &{}\\
{}           &{}                         &{}                    &{}}
    \end{xy}
    &
    \begin{xy}
  \xymatrix@!=10pt{%
{}\dd        &{}\dd                      &{}\dd                 &{}\dd\\
\lambda\dr \c&\frac{\lambda+1}{\mu}\dr \y&\frac{\nu}{\mu+1}\dr  &{}\d\\
\lambda\ddr  &\frac{\lambda}{\mu}\ddr    &\frac{\kappa}{\mu}\ann{\mbox{*}}\ddr&{}\dd\\
{}\dr        &{}\dr                      &{}\dr             \y\c&{}\d\\
{}\ddr       &{}\ddr                     &{}\ddr                &{}\dd\\
{}\dr      \y&{}\dr                    \c&{}\dr                 &{}\d\\
{}\ddr       &{}\ddr                     &{}\ddr                &{}\dd\\
{}\put{$k-1$}&{}\put{$k$}                &{}\put{$k+1$}         &{}\\
{}           &{}                         &{}                    &{}}
    \end{xy}\\
    a. Knuth relation of the first kind. &
    b. Knuth relation of the second kind.
  \end{tabular}
  \caption{Knuth equivalent columns}
    \label{fig:Knuth-columns}
  \end{figure}

  We first show that, in the situation of Figure~\ref{fig:Knuth-columns}.b,
  either $\mu=\lambda+1$ or $\nu=\mu+1$. To this end, suppose that
  $\mu\neq\lambda+1$. Since the fillings in the rectangles below and to the
  left of $*$ are equal, with the exception of the position of an empty column,
  $\mu=\kappa$. Applying forward rule F2, it follows that
  $\nu=\mu\cup\lambda+1=\mu+1$.

  It remains to show that the two growth diagrams differ in exactly one of the
  two partitions on the upper border.  If $\mu=\lambda+1$ then, by the forward
  rules, all partitions on the right side above the circle in column $k-1$ are
  identical, including the partitions on the upper border in this column.

  If, however, $\nu=\mu+1$, we have to use an inductive argument to show our
  claim. Consider the following part of a growth diagram, where both cells are
  assumed to be empty and $r\neq s$:

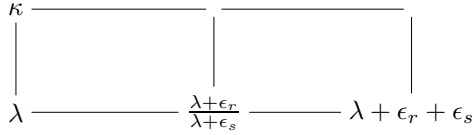
  \begin{equation*}
    \begin{xy}
      \xymatrix@!C=50pt{%
        \kappa\dr&
        {}\dr&
        {}\d\\
        \lambda\r&
        \frac{\lambda+\epsilon_r}{\lambda+\epsilon_s}\r&
        \lambda+\epsilon_r+\epsilon_s}
    \end{xy}
  \end{equation*}

  We will show that no matter what value $\kappa$ has, there will be again
  exactly one corner at the top where the partitions of the two growth diagrams
  differ:

  Suppose first, that $\kappa=\lambda+\epsilon_r$ and $s=r+1$. Applying the
  forward rules F2 and F3 we obtain
  
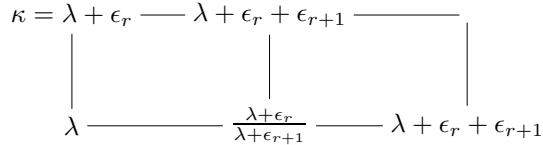
\begin{equation*}
    \begin{xy}
      \xymatrix@!C=50pt{%
        \kappa=\lambda+\epsilon_r\dr&
        \lambda+\epsilon_r+\epsilon_{r+1}\dr&
        {}\d\\
        \lambda\r&
        \frac{\lambda+\epsilon_r}{\lambda+\epsilon_{r+1}}\r&
        \lambda+\epsilon_r+\epsilon_{r+1}}
    \end{xy}
  \end{equation*}
  Thus the partitions labelling corners on the right of column $k-1$ of the two
  growth diagrams coincide from hereon, and the claim follows.

  If $\kappa=\lambda+\epsilon_r$, but $s\neq r+1$, we obtain (remember $s\neq r$)
  \begin{equation*}
    \begin{xy}
      \xymatrix@!C=50pt{%
        \kappa=\lambda+\epsilon_r\dr&
        \frac{\lambda+\epsilon_r+\epsilon_{r+1}}{\lambda+\epsilon_r+\epsilon_s}\dr&
        \lambda+\epsilon_r+\epsilon_{r+1}+\epsilon_s\d\\
        \lambda\r&
        \frac{\lambda+\epsilon_r}{\lambda+\epsilon_s}\r&
        \lambda+\epsilon_r+\epsilon_s}
    \end{xy}
  \end{equation*}
  
  Finally, if $\kappa=\lambda+\epsilon_t$, with all of $r$, $s$ and $t$ being
  different, we obtain
  \begin{equation*}
    \begin{xy}
      \xymatrix@!C=50pt{%
        \kappa=\lambda+\epsilon_t\dr&
        \frac{\lambda+\epsilon_r+\epsilon_t}{\lambda+\epsilon_s+\epsilon_t}\dr&
        \lambda+\epsilon_r+\epsilon_s+\epsilon_t\d\\
        \lambda\r&
        \frac{\lambda+\epsilon_r}{\lambda+\epsilon_s}\r&
        \lambda+\epsilon_r+\epsilon_s}
    \end{xy}
  \end{equation*}

  In both cases, the situation in the corners of the top row is the same as in
  the bottom row, and the claim follows.

  \textcircled{\footnotesize$\Leftarrow$} We assume now that the sequences of
  partitions along the upper border of the two Knuth equivalent growth diagrams
  differ in exactly one partition, namely between columns $k$ and $k+1$, and
  that either the three columns $k-1$, $k$ and $k+1$ or $k$, $k+1$ and $k+2$
  are shape equivalent to a triangle.

  Distinguishing between several cases, we consider the effect of applying the
  backward rules B1 to B4 to the partitions labelling the corners of columns
  $k-1$, $k$ and $k+1$. We will show that, in each row, if not all partitions
  coincide, either the partitions between columns $k$ and $k+1$ or those
  between column $k-1$ and $k$ differ, and all others coincide. Working our way
  from the top of the growth diagrams to the bottom, we can distinguish three
  stages: in Stage~1, the partitions between columns $k$ and $k+1$ differ. This
  stage can be followed either by Stage~2 or by Stage~3. In Stage~2, the
  partitions between columns $k-1$ and $k$ differ and only Stage~3 can follow.
  Finally, Stage~3 corresponds to those rows, that lie between the entries
  which are swapped by the Knuth transformation. It is followed by those rows
  in which all partitions coincide.

  Throughout the rest of this proof we will maintain $r<s$. To make it easier
  to keep track of the various situations and make the description more
  concise, we omit cases which may be obtained by exchanging the two growth
  diagrams. The backward rules we apply are printed into the corresponding
  cells wherever appropriate.

  \begin{enumerate}
  \item During this stage, the partitions of the two growth diagrams between
    columns $k$ and $k+1$ differ. The Cases~c and f lead directly to Stage~3,
    whereas Case~e leads to Stage~2. All other cases stay within Stage~1.
    \begin{equation*}
      \begin{xy}
        \xymatrix@!C=50pt{%
          \lambda-\epsilon_r-\epsilon_s\dr&
          \frac{\lambda-\epsilon_r}{\lambda-\epsilon_s}\dr&
          \lambda\d\\
          {}\r&{}\r&\kappa}
      \end{xy}
    \end{equation*}
    \begin{enumerate}
    \item $\kappa=\lambda-\epsilon_t$ and $r$, $s$ and $t$ are all different.

      \begin{equation*}
        \begin{xy}
          \xymatrix@!C=50pt{%
            \lambda-\epsilon_r-\epsilon_s\dr&
            \frac{\lambda-\epsilon_r}{\lambda-\epsilon_s}\dr\B2&
            \lambda\d\B2\\
            \kappa-\epsilon_r-\epsilon_s\r&
            \frac{\kappa-\epsilon_r}{\kappa-\epsilon_s}\r&
            \kappa=\lambda-\epsilon_t}
        \end{xy}
      \end{equation*}
      
    \item $\kappa=\lambda-\epsilon_r$ and $r>1$. Since $r-1<r<s$,

      \begin{equation*}
        \begin{xy}
          \xymatrix@!C=50pt{%
            \lambda-\epsilon_r-\epsilon_s\dr&
            \frac{\lambda-\epsilon_r}{\lambda-\epsilon_s}\dr\fB21&
            \lambda\d\fB32\\
            \kappa-\epsilon_{r-1}-\epsilon_s\r&
            \frac{\kappa-\epsilon_{r-1}}{\kappa-\epsilon_s}\r&
            \kappa=\lambda-\epsilon_r}
        \end{xy}
      \end{equation*}

    \item $\kappa=\lambda-\epsilon_r$ and $r=1$.

      \begin{equation*}
        \begin{xy}
          \xymatrix@!C=50pt{%
            \lambda-\epsilon_1-\epsilon_s\y\dr&
            \frac{\lambda-\epsilon_1}{\lambda-\epsilon_s}\dr\c\fB24&
            \lambda\d\fB42\\
            \kappa-\epsilon_s\r&
            \frac{\kappa}{\kappa-\epsilon_s}\r&
            \kappa=\lambda-\epsilon_1}
        \end{xy}
      \end{equation*}

    \item $\kappa=\lambda-\epsilon_s$ and $s\neq r+1$.

      \begin{equation*}
        \begin{xy}
          \xymatrix@!C=50pt{%
            \lambda-\epsilon_r-\epsilon_s\dr&
            \frac{\lambda-\epsilon_r}{\lambda-\epsilon_s}\dr\fB32&
            \lambda\d\fB23\\
            \kappa-\epsilon_r-\epsilon_{s-1}\r&
            \frac{\kappa-\epsilon_r}{\kappa-\epsilon_{s-1}}\r&
            \kappa=\lambda-\epsilon_s}
        \end{xy}
      \end{equation*}

      Note that, since $r<s$, we have automatically $s>1$.

    \item $\kappa=\lambda-\epsilon_s$ and $s=r+1$, $r>1$.

      In this case we also have to consider the partitions labelling the
      corners of columns $k-1$ and $k+2$.

      \begin{equation*}
        \begin{xy}
          \xymatrix@!C=50pt{%
            \lambda-\epsilon_r-\epsilon_{r+1}-\epsilon_t\dr&             
            \lambda-\epsilon_r-\epsilon_{r+1}\dr& 
            \frac{\lambda-\epsilon_r}{\lambda-\epsilon_{r+1}}\dr\B3&
            \lambda\dr\fB23&
            \lambda+\epsilon_u\d\\
            \mu\r&
            \frac{\kappa-2\epsilon_r}{\kappa-\epsilon_{r-1}-\epsilon_r}\r&
            \kappa-\epsilon_r\r&
            \kappa=\lambda-\epsilon_{r+1}\r&
            {}}
        \end{xy}
      \end{equation*}

      It turns out that in this case columns $k$, $k+1$ and $k+2$ cannot be
      shape equivalent to a triangle in both fillings.  We first observe that
      this would force $u=r+1$: if $u<r+1$, columns $k$, $k+1$ and $k+2$ form a
      row in the diagram corresponding to circles (i.e., the partitions on
      top), since the partitions $\lambda-\epsilon_r-\epsilon_{r+1}$,
      $\lambda-\epsilon_r$, $\lambda$ and $\lambda+\epsilon_u$ differ in
      \emph{columns} with strictly increasing indices.  Similarly, if $u>r+1$,
      columns $k$, $k+1$ and $k+2$ form a column in the diagram corresponding
      to crosses, because then the partitions
      $\lambda-\epsilon_r-\epsilon_{r+1}$, $\lambda-\epsilon_{r+1}$, $\lambda$
      and $\lambda+\epsilon_u$ differ in \emph{rows} with strictly increasing
      indices.

      However, if $u=r+1$, we have a cell
      \begin{equation*}
        \begin{xy}
          \xymatrix{%
            \lambda    \dr&\lambda+\epsilon_{r+1}\d\\
            \lambda-\epsilon_{r+1}\r &{}}
        \end{xy}
      \end{equation*}
      which is impossible, considering the local rules.

      We conclude that columns $k-1$, $k$ and $k+1$ must be shape equivalent to
      a triangle. By reasoning similar to above we find $t=r$ and therefore
      $\mu=\kappa-\epsilon_{r-1}-2\epsilon_r$ by B3 and B2 respectively.

    \item $\kappa=\lambda-\epsilon_s$ and $s=r+1$, $r=1$.

      This case is very similar to the preceding one.  Again we find that
      $t=r$, and conclude that the cell in the top row and column $k-1$ will
      now contain a circle, the cell in column $k$ a cross,
      $\mu=\kappa-2\epsilon_1$:

      
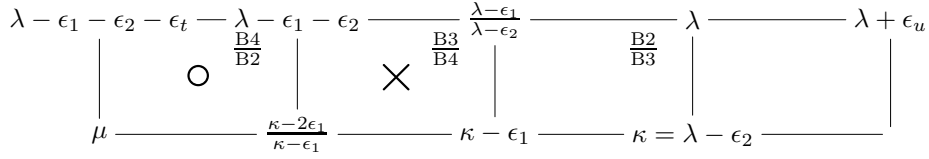
\begin{equation*}
        \begin{xy}
          \xymatrix@!C=50pt{%
            \lambda-\epsilon_1-\epsilon_2-\epsilon_t\dr\c&
            \lambda-\epsilon_1-\epsilon_2\dr\y\fB42& 
            \frac{\lambda-\epsilon_1}{\lambda-\epsilon_2}\dr\fB34&
            \lambda\dr\fB23&
            \lambda+\epsilon_u\d\\
            \mu\r&
            \frac{\kappa-2\epsilon_1}{\kappa-\epsilon_1}\r&
            \kappa-\epsilon_1\r&
            \kappa=\lambda-\epsilon_2\r&
            {}}
        \end{xy}
      \end{equation*}

    \end{enumerate}

  \item Within this stage, the partitions between columns $k-1$ and $k$ of the
    two growth diagrams differ.  We show that cases in this stage can only be
    followed by cases of Stage~2 or 3.
    
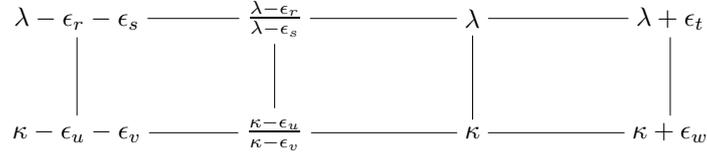
\begin{equation*}
      \begin{xy}
        \xymatrix@!C=50pt{%
          \lambda-\epsilon_r-\epsilon_s\dr& 
          \frac{\lambda-\epsilon_r}{\lambda-\epsilon_s}\dr&
          \lambda\dr&
          \lambda+\epsilon_t\d\\
          \kappa-\epsilon_u-\epsilon_v\r& 
          \frac{\kappa-\epsilon_u}{\kappa-\epsilon_v}\r&
          \kappa\r&
          \kappa+\epsilon_w}
        \end{xy}
    \end{equation*}
    This situation occurs after Case~1e and thus can only occur if columns
    $k-1$, $k$ and $k+1$ are shape equivalent to a triangle in both fillings.
    Throughout this subcase, in the situation of the diagram above, we maintain
    $r<t\le s$.  Note that this condition is satisfied by the partitions along
    the bottom of the diagram in Case~1e.  We are going to show that we also
    have $u<w\le v$ (in Cases~a,b, and c) or that we continue with one of the
    cases of Stage~3.
    \begin{enumerate}
    \item $\kappa=\lambda-\epsilon_x$ and $r$, $s$ and $x$ are all different.

      By B2 we obtain
      $\kappa-\epsilon_u=\lambda-\epsilon_r\cap\lambda-\epsilon_x$ and thus
      $u=r$. Similarly, we have $v=s$. To determine $w$, we consider the cell
      
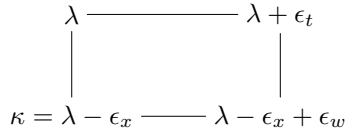
\begin{equation*}
        \begin{xy}
          \xymatrix{%
            \lambda    \dr&\lambda+\epsilon_t\d\\
            \kappa=\lambda-\epsilon_x\r &\lambda-\epsilon_x+\epsilon_w}
        \end{xy}
      \end{equation*}
      If $w=x$, we obtain (remember $t>1$) by B3 $x=t-1$.  Since $r\neq x=w$
      and $r\le t-1=w$ we also have $u=r<w=t-1\le s=v$.

      If, however, $w\neq x$, then the local rules force $w=t$.

    \item $\kappa=\lambda-\epsilon_s$.

      By B2 we obtain $u=r$, by B3 that $v=s-1$. Now consider
      
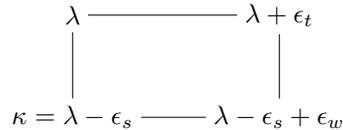
\begin{equation*}
        \begin{xy}
          \xymatrix{%
            \lambda    \dr&\lambda+\epsilon_t\d\\
            \kappa=\lambda-\epsilon_s\r &\lambda-\epsilon_s+\epsilon_w}
        \end{xy}
      \end{equation*}
      If $w=s$, the local rules imply $t=s+1$, which contradicts our assumption
      $t\le s$. Thus, $w\neq s$, which entails $w=t$.  Since $t=s$ is
      impossible by the local rules, we have indeed $u=r<w=t\le s-1=v$.

    \item $\kappa=\lambda-\epsilon_r$ and $r>1$.

      By B3 we obtain $u=r-1$ and by B2 $v=s$. To determine $w$, we consider the cell
      \begin{equation*}
        \begin{xy}
          \xymatrix{%
            \lambda    \dr&\lambda+\epsilon_t\d\\
            \kappa=\lambda-\epsilon_r\r &\lambda-\epsilon_r+\epsilon_w}
        \end{xy}
      \end{equation*}
      If $w=r$ we immediately have $u=r-1<w=r\le v=s$. If $w\neq x$, we obtain
      $w=t$ and $u=r-1<w=t\le s=v$.

    \item $\kappa=\lambda-\epsilon_r$ and $r=1$.

      In this case the growth diagram looks as follows:
      \begin{equation*}
        \begin{xy}
          \xymatrix@!C=50pt{%
            \lambda-\epsilon_1-\epsilon_s\dr\y& 
            \frac{\lambda-\epsilon_1}{\lambda-\epsilon_s}\dr\fB24\c&
            \lambda\dr\fB42&
            \lambda+\epsilon_t\d\\
            \lambda-\epsilon_1-\epsilon_s\r& 
            \frac{\lambda-\epsilon_1}{\lambda-\epsilon_1-\epsilon_s}\r&
            \kappa=\lambda-\epsilon_1\r&
            {}}
        \end{xy}
      \end{equation*}
      We are thus lead to a situation of Stage~3.
    \end{enumerate}

  \item In this final stage, the growth diagram is of one of the following
    forms:
    \begin{enumerate}
    \item $u\neq r$
      \begin{equation*}
        \begin{xy}
          \xymatrix@!C=50pt{%
            \lambda-\epsilon_r\dr& 
            \frac{\lambda}{\lambda-\epsilon_r}\dr&
            \lambda\d\\
            \lambda-\epsilon_u-\epsilon_r\r& 
            \frac{\lambda-\epsilon_u}{\lambda-\epsilon_u-\epsilon_r}\r&
            \lambda-\epsilon_u}
        \end{xy}
      \end{equation*}
    \item $r>1$
      \begin{equation*}
        \begin{xy}
          \xymatrix@!C=50pt{%
            \lambda-\epsilon_r\dr& 
            \frac{\lambda}{\lambda-\epsilon_r}\dr&
            \lambda\d\\
            \lambda-\epsilon_{r-1}-\epsilon_r\r& 
            \frac{\lambda-\epsilon_r}{\lambda-\epsilon_{r-1}-\epsilon_r}\r&
            \lambda-\epsilon_r}
        \end{xy}
      \end{equation*}
    \item 
      \begin{equation*}
        \begin{xy}
          \xymatrix@!C=50pt{%
            \lambda-\epsilon_1\dr\c& 
            \frac{\lambda}{\lambda-\epsilon_1}\dr\y&
            \lambda\d\\
            \lambda-\epsilon_1\r& 
            \lambda-\epsilon_1\r&
            \lambda-\epsilon_1}
        \end{xy}
      \end{equation*}
    \end{enumerate}
  \end{enumerate}

  By considering the different cases as detailed above, the reader can convince
  herself of the validity of the claim.
\end{proof}

\begin{lem}\label{lem:jdt-difference}
  Let $R$ and $S$ be a partial Young tableaux obtained from each other by
  interchanging the entries $k$ and $k+1$.  Suppose that the entries $k-1$, $k$
  and $k+1$ of both $R$ and $S$ are shape equivalent to a triangle.  Then
  $jdt(S)$ can be obtained from $jdt(R)$ by interchanging entries $k-1$ and
  $k$.

  Similarly, suppose that $S$ and $R$ are obtained from each other by
  interchanging the entries $k-1$ and $k$, and entries $k-1$, $k$ and $k+1$ of
  both $R$ and $S$ are shape equivalent to a triangle.  Then $jdt(S)$ can be
  obtained from $jdt(R)$ by either interchanging entries $k-2$ and $k-1$ or
  entries $k-1$ and $k$.
\end{lem}
\begin{proof}
  Consider the effect of applying the algorithm $jdt$ to $R$, using the
  description via slides on tableaux.  If the two cells whose entries differ in
  $R$ and $S$ are not adjacent, their entries will never be compared by the
  algorithm.  The statement then follows since $jdt$ subtracts one from all
  entries of the tableau.

  Now consider the case that $k$ and $k+1$ are adjacent.  We first show that
  the cells containing $k-1$, $k$ and $k+1$ must be arranged as follows in $R$
  and $S$:
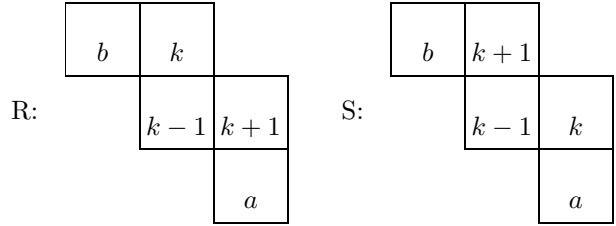
  \begin{equation*}
    \Yboxdim28pt
    \def\km{k-1}
    \def\kp{k+1}
    \text{R:}\quad\young(bk,:\km\kp,::a)\quad\quad
    \text{S:}\quad\young(b\kp,:\km k,::a)
  \end{equation*}
  If $k-1$ were in a higher row, the entries $1$, $2$ and $3$ would form a
  single row in $jdt^{k-2}(R)$.  But this is incompatible with the hypothesis
  that they are shape equivalent to a triangle.  Similarly, if $k-1$ were in a
  lower row, $1$, $2$ and $3$ would form a column in $jdt^{k-2}(S)$, which
  again contradicts the hypothesis.

  We now observe that both $a$ and $b$ must be smaller than $k-1$, since
  entries in rows and columns are strictly increasing.  Thus, in this case,
  $jdt$ can never swap the empty cell with the cell containing $k-1$.  Taking
  into account that $jdt$ subtracts one from all entries, we conclude that
  $jdt(S)$ can be obtained from $jdt(R)$ by exchanging entries $k-1$ and $k$.

  We now consider the case that entries $k-1$ and $k$ were swapped to obtain
  $S$ from $R$, and the cells containing $k-1$ and $k$ are adjacent.  Then, the
  cells containing $k-1$, $k$ and $k+1$ must be arranged as follows:
  \begin{equation*}
    \Yboxdim28pt
    \def\km{k-1}
    \def\kp{k+1}
    \text{R:}\quad\young(b,\km\kp,*ka)\quad\quad
    \text{S:}\quad\young(b,k\kp,*\km a)
  \end{equation*}
  In this situation, both $a$ and $b$ must be greater than $k+1$.  Thus, if
  during $jdt$ the empty cell is moved to the place marked with a $*$ above,
  after two steps of the algorithm we have the following situations:
  \begin{equation*}
    \Yboxdim28pt
    \def\km{k-1}
    \def\kp{k+1}
    \quad\young(b,\kp *,\km ka)\quad\quad
    \quad\young(b,k*,\km\kp a)
  \end{equation*}
  Therefore, only entries $k$ and $k+1$ are exchanged.  Taking into account
  that $jdt$ subtracts one from all entries, we conclude that $jdt(S)$ can be
  obtained from $jdt(R)$ by exchanging entries $k-1$ and $k$.
\end{proof}

Let us restate Proposition~\ref{prop:equivalence} for easier reference:
\begin{prop}
  Consider the following two fillings:
  \begin{gather*}
    \young(\alpha\beta\delta)
    \quad\text{and}\quad
    \young(\alpha\gamma\delta)
  \end{gather*}
  and suppose furthermore that $\beta$ and $\gamma$ are Knuth equivalent.
  Then, applying $\J$ to both fillings we obtain
  \begin{gather*}
    \young(\alphaP\betaP\deltaP)
    \quad\text{and}\quad
    \young(\alphaP\gammaP\deltaP)
  \end{gather*}
  where $\alpha^\prime$ has exactly one column less than $\alpha$ and
  $\delta^\prime$ has exactly one more column than $\delta$.  In this
  situation, $\beta^\prime$ and $\gamma^\prime$ are Knuth equivalent.

  Similarly, consider
  \begin{gather*}
    \young(\delta,\beta,\alpha)
    \quad\text{and}\quad
    \young(\delta,\gamma,\alpha)
  \end{gather*}
  and suppose furthermore that $\beta$ and $\gamma$ are dual Knuth equivalent.
  Then, applying $\J$ to both fillings we obtain
  \begin{gather*}
    \young(\deltaP,\betaP,\alphaP)
    \quad\text{and}\quad
    \young(\deltaP,\gammaP,\alphaP)
  \end{gather*}
  where $\alpha^\prime$, $\beta^\prime$ and $\delta^\prime$ have as many rows
  as $\alpha$, $\beta$ and $\delta$ respectively.  In this situation,
  $\beta^\prime$ and $\gamma^\prime$ are dual Knuth equivalent.
\end{prop}

\begin{proof}
  We commence by proving the case of standard fillings.  To prove the first
  statement, we observe that $\beta$ and $\gamma$ can be transformed one into
  the other using a sequence of Knuth transformations.  By transitivity, we
  only need to consider the case where $\beta$ and $\gamma$ differ by a single
  Knuth relation.  In this situation we can apply
  Proposition~\ref{prop:difference} and Lemma~\ref{lem:jdt-difference}: by
  Proposition~\ref{prop:difference}, the upper borders of $\beta$ and $\gamma$
  differ in a single partition, say between columns $k$ and $k+1$.  Therefore,
  the corresponding partial Young tableaux can be obtained from each other by
  interchanging $k$ and $k+1$ and we can apply Lemma~\ref{lem:jdt-difference}.
  It follows that the partitions labelling the upper border of $\beta'$ and
  $\gamma'$ also differ in a single position.  Thus, we can apply the converse
  direction of Proposition~\ref{prop:difference} to conclude that the fillings
  $\young(\alphaP\betaP\deltaP)$ and $\young(\alphaP\gammaP\deltaP)$ differ
  only by a single Knuth relation.

  The second statement follows by similar reasoning directly from
  Proposition~\ref{prop:difference}.  We only have to observe that the notions
  of Knuth equivalence and dual Knuth equivalence are connected by reflecting
  diagrams about the main diagonal.

  \Yboxdim13pt It remains to reduce the case of arbitrary fillings to the
  standard case. Let $\young(\alphaB\betaB\deltaB)$ be the partial
  standardisation of the filling $\young(\alpha\beta\delta)$ described in
  Lemma~\ref{lem:commutation}, i.e., all columns but the first of $\alpha$ are
  standardised, similarly for $\young(\alphaB\gammaB\deltaB)$.  As before, we
  can now assume that $\betaB$ and $\gammaB$ differ by a single Knuth relation,
  and deduce that $\young(\alphaBP\betaBP\deltaBP)$ and
  $\young(\alphaBP\gammaBP\deltaBP)$ differ only by a single Knuth relation.
  By Lemma~\ref{lem:commutation} we have
  $std\left(\J\left(\young(\alpha\beta\delta)\right)\right)=
  std\left(\young(\alphaBP\betaBP\deltaBP)\right)$ and
  $std\left(\J\left(\young(\alpha\gamma\delta)\right)\right)=
  std\left(\young(\alphaBP\gammaBP\deltaBP)\right)$, and the claim follows.
\end{proof}


\providecommand{\bysame}{\leavevmode\hbox to3em{\hrulefill}\thinspace}
\providecommand{\href}[2]{#2}

\end{document}